\documentclass[12pt]{article}

\usepackage{amsmath,amssymb,amsthm}
\usepackage{xcolor}
\usepackage{hyperref}
\usepackage{bm}
\usepackage{graphicx}
\usepackage{subcaption}
\usepackage{url}
\newtheorem{theorem}{Theorem}[section]    
\newtheorem{corollary}[theorem]{Corollary}
\newtheorem{proposition}[theorem]{Proposition}
\newtheorem{definition}[theorem]{Definition}

\theoremstyle{remark}
\newtheorem{remark}[theorem]{Remark}

\title{Spherical Poisson Needlets with Shrinking Bandwidth}

\author{Mattia Castaldo \and Claudio Durastanti \\
	Department of Basic Science and Applied to Engineering, \\
	Sapienza University of Rome, Italy \\
	\texttt{claudio.durastanti@uniroma1.it}
}

\date{} 

\begin{document}
	
	\maketitle
\begin{abstract}
Flexible bandwidth needlets provide a localized multiscale framework with scale-adaptive frequency resolution, enabling effective analysis of spherical Poisson random fields exhibiting spatial inhomogeneity and scale variation. We establish here quantitative Central Limit Theorems for finite-dimensional distributions of spherical Poisson needlets and for the related Poisson needlet coefficients constructed via needlets with shrinking bandwidth on the sphere, and using Stein–Malliavin techniques, we derive explicit rates of normal approximation.
In addition, we study the functional convergence of the associated needlet-based random fields. Indeed, our framework provides quantitative control on the limiting behavior in appropriate function spaces. Together, these results offer rigorous probabilistic guarantees for high-resolution spherical data modeling under Poisson sampling.
\end{abstract}
\section{Introduction}
Needlets on the sphere are wavelets characterized by strong concentration properties in both the spatial and frequency domains. Standard needlets, introduced in \cite{npw1,npw2}, fix a band factor $B>1$ so that, at resolution level $j \geq 1$, each needlet $\psi_{j,k}$ is essentially supported on a spherical cap of radius $O(B^{-j})$ and bandlimited to harmonics in the narrow interval $[B^{j-1},B^{j+1}]$, see also \cite{gm2}.  The scaling $B^{j}$ thus serves as a bandwidth, governing both the scale of spatial localization and the width of the frequency window. Related wavelet constructions have also been introduced in the literature, as for instance \cite{bgsw24}, where exact cubature rules are replaced with quasi‐Monte Carlo designs while preserving Sobolev convergence rates. Another example are the various directional needlet frameworks (e.g. \cite{mcewen,schoppert}), where the wavelets are additionally indexed by orientation so as to achieve simultaneous localization in scale, position, and direction on the sphere. 

Standard needlets have been used in many statistical applications, such as nonparametric regression or density estimation on $\mathbb{S}^2$, where one shrinks or thresholds these needlet coefficients to adapt to unknown smoothness, achieving near–minimax risk bounds (see \cite{bkmpAoSb,dur2,monnier}), with some extension to nonparametric random‐coefficients binary‐choice model on the sphere (see \cite{gp18}).
In \cite{lmejs}, the authors introduce a needlet‐based estimator of the angular bispectrum on the sphere, proving its asymptotic normality and showing how it provides a highly localized, computationally efficient tool for detecting and quantifying non-Gaussianity with explicit error bounds. In \cite{dlmejs}, a semiparametric needlet–Whittle estimator has been developed to exploit the tight localization of needlets to yield efficient, nearly optimal estimates of the angular power spectrum under minimal smoothness assumptions. \\
Since the foundational work of \cite{peczen}, Stein–Malliavin techniques have undergone continuous theoretical development in the Poisson setting, with significant advancements presented in works such as \cite{dvz}, and more recent contributions in \cite{be24,tt}.
In the Poisson framework, needlets have been extensively applied across various applications. Indeed, these same localization and frame properties have been harnessed, via Malliavin–Stein techniques, to prove quantitative Gaussian approximations for both linear and nonlinear needlet coefficients built over spherical Poisson fields, thus providing high‐dimensional central limit theorems and explicit convergence rates for thresholding and inference in the regime of growing count intensity. In \cite{dmp14} suitably normalized linear needlet sums converge to the normal law with explicit error bounds; in \cite{bdmp}  a broad class of nonlinear functionals have been discussed, proving quantitative central limit theorems for them. Finally, in \cite{bd}, the authors showed that high-frequency U-statistics whose kernel functions live in classical nonparametric spaces such as Besov classes converge in distribution to a Gaussian law. This result not only bridges the gap between geometric multiscale analysis and traditional nonparametric inference—by tying the smoothness of the kernel in Besov norms directly to convergence rates—but also provides a concrete toolbox for constructing asymptotically normal estimators and test statistics in settings ranging from density and intensity estimation to hypothesis testing on curved domains.Together, these results establish a powerful multiscale toolkit for statistical inference on point-process data over curved domains, marrying the localization and frame properties of needlets with modern Malliavin–Stein techniques. \\
A completely different approach to analyzing Poisson waves on the sphere was developed in \cite{bdmt24}, where the authors introduced a model of Poisson-driven random waves by expanding the underlying point process in the spherical harmonic basis. Quantitative central limit theorems were then derived for finite-dimensional harmonic distributions and functional convergence, capturing the interplay between the Poisson intensity and the eigenvalue growth. In this work, we bridge the gap between this harmonic expansion framework and the flexible bandwidth needlet construction, combining the strengths of both approaches for refined multiscale analysis of Poisson spherical fields.\\
Flexible bandwidth needlets extend the standard needlet construction by replacing the fixed exponential scaling $B^j$, for some base $B > 1$, with a more general, adaptable scaling sequence $\{S_j: j \geq 1\}$. While standard needlets localize frequencies around bands growing geometrically as powers of $B$, flexible bandwidth needlets allow the frequency bands to increase at rates tailored to the data or application, potentially slower or faster than a pure exponential. This flexibility enables finer control over the localization properties in both frequency and space, making the framework well suited for handling inhomogeneous or scale-varying phenomena on the sphere, as detailed in \cite{dmt24}.
In \cite{d25}, it has been shown that a proper choice of $S_j$ to enforce frequency concentration (the so called shrinking regime) is given by
\[
S_j = \gamma(j)/j^{p},
\] 
where $p\in[0,1]$ reflects regularity and $\gamma(j)$ s a slowly varying function in the sense of Karamata (see \cite{bingham}), meaning that for all  $\tau>0$
\[
\underset{j \rightarrow \infty}{\lim}\frac{\gamma(\tau j)}{\gamma( j)}=1.
\]
Indeed, these shrinking needlets retain the same tight‐frame and near‐exponential localization properties, yet allow the resolution to adapt dynamically to data density or higher‐order smoothness. Fixed a resolution level $j$, the needlet projector (or kernel) is given by
\[
\Phi_j(x,y)=\sum_{\ell = S_{j-1}}^{S_{j+1}} b_j^2 \left(\ell\right) Z_{\ell} (\langle x, y\rangle),
\]
while the shrinking needlets are defined by:
\[
\psi_{j,k}= \sqrt{\lambda_{j,k} }\sum_{\ell = S_{j-1}}^{S_{j+1}} b_j \left(\ell\right) Z_{\ell} (\langle x, \xi_{j,k}\rangle),
\]
where $\xi_{j,k}$ and $\lambda_{j,k}$ are properly chosen cubature points and weights, $b_j: \mathbb{R}^+ \mapsto \mathbb{R}^+$ is a weight function and $Z_{\ell}$ is the standard harmonic projector onto the multipole $\ell$. 
%
As aforementioned, in this paper we discuss the construction of spherical Poisson needlets, built by projecting a random Poisson measure onto a shrinking-needlet frame, whose bandwidth contracts with increasing frequency. To define spherical Poisson needlets, we convey the notion of random phases on the sphere by using the model of spherical Poisson waves introduced by \cite{bdmt24}, and defined as superpositions of spherical harmonics centered at random locations distributed according to a Poisson process on $\mathbb{S}^2$. Here we can introduce this construction as follows. For a fixed level $j \geq 1$, the spherical Poisson needlet $\Psi_{j;t}(x)$,  is defined by
\[
\Psi_{j;t}(x)= \frac{1}{\sqrt{\nu_t}\sigma_j}\sum_{i=1}^{N_t(\mathbb{S}^2)} \Phi_{j}(x,x_i), \quad x \in \mathbb{S}^2,
\]
where $\lbrace x_1, x_2,\ldots \rbrace$ are independent random locations on the sphere whose amount is random and sampled from a Poisson point process on the sphere $N_t\left(\mathbb{S}^2\right)$ with intensity $\nu_t$, and $\sigma_j$ is a normalization factor. This formulation captures the interplay between concentration of the needlet construction and randomness via the integration of deterministic needlet projector (that is, a weighted sum of Laplace spherical eigenfunctions) against a Poisson random measure on the sphere.
Shrinking needlets occupy an intermediate position between global spherical harmonics and fixed-bandwidth needlets: they inherit the exact spectral localization of harmonics while preserving the strong spatial concentration and tight‐frame properties of needlets. This intermediate construction gives us finer control over frequency selection, allowing the bandwidth to narrow as the resolution level grows and thereby delivering sharper asymptotic approximations for both individual coefficients and their joint behavior. As a result, we can track the interplay between Poisson intensity, spectral decay, and spatial localization more precisely than in the classical needlet or harmonic settings, yielding quantitative central‐limit theorems and functional convergence results with explicitly controlled error terms.\\
Our analysis begins with a detailed study of the convergence in law of the needlet Poisson coefficients themselves, establishing quantitative central limit theorems that describe their joint asymptotic behavior. Building on this foundation, we then examine the finite-dimensional distributions of shrinking needlets evaluated at a fixed collection of points on the sphere. For the univariate case, we prove a central limit theorem with explicit convergence rates linked to the growth of scaling and the Poisson intensity. This extends naturally to the multivariate setting, where the dimension grows with the observation scale, and quantitative bounds ensure asymptotic Gaussianity under suitable growth conditions. Crucially, we also provide a rigorous link between the convergence of the coefficients and that of the finite-dimensional distributions, highlighting how the two perspectives complement each other. Finally, we investigate functional convergence results in $L^2(\mathbb{S}^2)$ and stronger Sobolev spaces, revealing subtle distinctions in topological effects and convergence rates.
\subsubsection*{Plan of the paper}
In Section \ref{sec:prel}, we introduce spherical Poisson waves and flexible bandwidth needlets, along with the notation and key theoretical results from \cite{dvz} that are essential for establishing our quantitative bounds. In Section \ref{sec:main},  
first we derive and discuss quantitative convergence results specifically for the shrinking needlet coefficients associated with the spherical Poisson wave model, 
then we establish quantitative central limit theorems for the finite‑dimensional vectors of spherical Poisson needlets evaluated at a fixed set of well‑separated locations on the sphere. Section \ref{sec:proofs} collects all the proofs.

\section{Preliminary results}\label{sec:prel}
In this section, we lay the groundwork for our study by introducing the framework of spherical Poisson waves and the associated needlet constructions with flexible bandwidth. We review the essential definitions, notation, and properties of these objects, drawing on key results from \cite{dvz} that will serve as the foundation for the aymptotic results developed later. \\
From now on, we adopt the following notational conventions throughout the article. Consider two sequences $\lbrace a_j:j \geq 1\rbrace$ and $\lbrace b_j:j \geq 1\rbrace$. When we write $a_j \sim b_j$, then the two sequences are asymptotically equivalent, that is, there exists a real -valued constant $C\neq0$ such that
\[
\underset{j \rightarrow \infty}{\lim}\frac{a_j}{b_j} \rightarrow C.
\]
Similarly, we use the notation $a_j \lesssim b_j$ to indicate that $ a_j$ is asymptotically bounded above by $b_j$, up to a constant factor independent of $j$.
 
\subsection{The construction of flexible bandiwidth needlets on the sphere}
Before introducing flexible bandwidth needlets, we briefly review some standard background material on harmonic analysis on the sphere. For further discussion and comprehensive details, we refer the reader to \cite{atki,MaPeCUP} and the references therein.
Let $x=\left(\vartheta,\varphi\right)\in\mathbb{S}^2$ be a location on the sphere; $\vartheta \in \left(0,\pi\right)$ is the colatitude, while $\mathcal{L}( \operatorname{d}x) = \sin \vartheta \operatorname{d}\vartheta \operatorname{d}\varphi$ is the spherical Lebesgue measure. To simplify notation, we will henceforth use $\mathcal{L} ( \operatorname{d}x)$ and the abbreviated form $\operatorname{d}x$ interchangeably, whenever no ambiguity arises. The set of spherical harmonics $\lbrace  Y_{\ell ,m} : \ell \geq 0; m=-\ell,\ldots,\ell\rbrace$ is an orthonormal basis for the space of square-integrable functions with respect to $\operatorname{d}x$ by $L^{2}\left( \mathbb{S}^{2}\right) =L^{2}\left(\mathbb{S}^{2},\mathcal{L}(\operatorname{d}x)\right)$; $\ell\in \mathbb{N}_0$ is the multipole number indicating the total degree of the homogeneous harmonic polynomial exactly reconstructed by $\lbrace Y_{\ell,m}: m=-\ell,\ldots,\ell\rbrace$ when restricted to the sphere, and $m=-\ell,\ldots, \ell$ is the azimuthal number, denoting the degeneracy of linearly independent spherical harmonics of degree $\ell$ distinguished by their angular dependence. Also, we can define
\[
\mathcal{K}_{\ell} = \bigoplus_{\ell^\prime = 0}^{\ell} \operatorname{spam}\lbrace Y_{\ell^\prime,m}:m=-\ell^{\prime},\ldots,\ell^{\prime}\rbrace,
\]
the linear space containing the restriction to the sphere of the
polynomials with degree at most $\ell$. It is possible to show (see for example \cite{atki}) that any square-integrable function $f\in L^{2}\left( \mathbb{S}%
^{2}\right) $ admits the harmonic expansion
\begin{equation*}
	f\left( x\right) =\sum_{\ell \geq 0}\sum_{m=-\ell}^{\ell}a_{\ell
		,m}Y_{\ell ,m}\left( x\right) ,  \label{eq:exp}
\end{equation*}%
where, for $\ell \geq 0$ and $m=-\ell,\ldots ,\ell$,
\begin{equation*}
	a_{\ell ,m}=\int_{\mathbb{S}^{d}}\overline{Y}_{\ell ,m}\left( x\right)
	f\left( x\right) \operatorname{d}x \in \mathbb{C},
\end{equation*}%
are the spherical harmonic coefficients. Spherical harmonic enjoy a very useful addition formula which we can use to define the projector $Z_\ell$ onto the multipole level $\ell$:
\begin{eqnarray*}
	Z_{\ell }\left( x_{1},x_{2}\right) &=&\sum_{m=-\ell}^{\ell}\overline{Y}%
	_{\ell ,m}\left( x_{1}\right) Y_{\ell ,m}\left( x_{2}\right)  \notag \\
	&=&\frac{2\ell +1}{4\pi} P_{\ell }\left( \left\langle x_{1},x_{2}\right\rangle \right) ,\quad
	\text{for }x_{1},x_{2}\in \mathbb{S}^{2},  \label{eq:harmproj}
\end{eqnarray*}%
where $\left\langle \cdot ,\cdot \right\rangle $ is the standard scalar
product over $\mathbb{R}^{3}$, and $P_{\ell }$ is
the Legendre polynomial of degree $\ell $. 
Spherical harmonic projectors feature the following \emph{reproducing kernel property}:
\begin{equation*}
	\int_{\mathbb{S}^{d}}Z_{\ell }\left( \langle x,y\rangle \right) Z_{\ell
		^{\prime }}\left( \langle y,z\rangle \right) \operatorname{d} y=Z_{\ell }\left( \langle
	x,z\rangle \right) \delta _{\ell ^{\prime }}^{\ell },  \label{eq:selfreprod}
\end{equation*}%
where $\delta _{\cdot }^{\cdot }$ is the Kronecker delta.
Note that projection of  any $%
f\in L^{2}\left( \mathbb{S}^{2}\right) $ over the space of multipole $\ell$ is given by
\begin{equation*}
	f_{\ell }\left( x\right) =Z_{\ell}\left[ f\right] \left( x\right) =\int_{%
		\mathbb{S}^{2}}Z_{\ell}\left( x,y \right) f\left( y\right)
	\operatorname{d} y=\sum_{m=-\ell}^{\ell}a_{\ell ,m}Y_{\ell ,m}\left( x\right) .
	\label{eq:proj}
\end{equation*}
While classical spherical harmonic analysis decomposes a function into globally supported, fixed‐bandwidth basis functions, flexible‐bandwidth needlets extend this framework by employing localized, multiscale filters that adaptively concentrate spectral energy within chosen frequency bands, enabling both precise spatial localization and tunable spectral resolution
The flexible bandwidth needlet kernel, as detailed by \cite{dmt24,d25}, can be
defined as a weighted sum of harmonic projector: for any $j=1,2,...$%
\begin{equation*}
	\Phi _{j}\left( x,y\right) 
	=\sum_{\ell \geq
		0}b^2_j\left( \ell\right) Z_{\ell }\left( \left\langle
	x,y\right\rangle \right) \text{ ,}
\end{equation*}%
where 
$\lbrace b_j(\cdot):j\geq 1 \rbrace$ is sequence of a weight functions $ b_j:\mathbb{R}^+%
\rightarrow \mathbb{R}^+$. Each flexible bandwidth system is characterized by a scale (or center) sequence $\left\{ S_{j}\right\}$, a strictly increasing sequence of positive real numbers, which identify also the center of the weight function support. More in detail,
for each $j \geq 1$, each weight function must satisfy three
properties:
\begin{enumerate}
	\item \emph{Compact support and normalization.} $b_j$ is compactly supported in $[S_{j-1},S_{j+1}].$ Moreove it holds that $b_{j}\left(S_{j-1}\right) = b_{j}\left(S_{j+1}\right) = 0$, and $b_j\left(S_j\right)=1$. 

	\item \emph{Differentiability.} $b_j$ is $C^{\infty }$; also it holds that 
\begin{equation*}
	\left \vert b_j^{(n)} (u) \right \vert \leq \frac{C_n}{\left(S_j-S_{j-1}\right)^n},
\end{equation*}
where $C_n>0$ and $b_j^{(n)}$ is the $n$-th derivative of $b_j$
 
 \item \emph{Partition of unity property.} For all $\ell \geq 1$, 
$\sum_{j}b_j^{2}(\frac{\ell })\equiv 1$.
\end{enumerate}
The function $b_j$ is strictly localized in frequency around the central value $S_j$, with support $\left[S_{j-1},S_{j+1}\right]$. This property is crucial for preserving the semi-orthogonality characteristic of needlets: in particular, the supports of the weight functions $b_j$ and $b_{j^\prime}$ are disjoint whenever $\left \vert j-j^{\prime}\right\vert \geq 2$.\\
The normalization condition $b_j(S_j) = 1$ ensures full preservation of the central frequency at scale $j$, enabling multiscale analysis through localization in narrow frequency bands centered at $S_j$. The smoothness of the infinitely differentiable functions $b_j$ guarantees spatial localization of the needlets, while bounds on their derivatives, tied to the bandwidth $[S_{j-1}, S_{j+1}]$, control the sharpness of transitions, balancing frequency and spatial localization. Finally, the partition of unity property ensures complete and stable coverage of the frequency domain, allowing faithful decomposition and reconstruction in $L^2(\mathbb{S}^2)$. \\
For more details, the reader is referred to \cite{dmt24} for a numerical scheme for constructing $b_j$, which serves as a smooth spectral bump function localized in the frequency band $\left[S_{j-1},S_{j+1}\right]$, closely related to the method proposed in \cite{bkmpAoS} for standard needlets and also tailored to form a smooth partition of unity in the frequency domain.

To rigorously formulate our results, we adopt the notation introduced by \cite{d25}, which is specifically designed to derive explicit bounds on concentration properties. For $j\geq 1 $, we define the \emph{bandwidth dilation factor}
\begin{equation*}
	 h_{j}:=\frac{S_{j+1}}{S_{j}},
\end{equation*}
which controls the relative enlargement of the centers across scales in the needlet system.  
In this work, we focus on the so-called \emph{shrinking needlets}, characterized by the property that the ratio between the bandwidth width $\left[S_{j-1}, S_{j+1}\right]$ and its central frequency $S_j$ tends to zero as $j \to \infty$. This implies increasingly sharp frequency localization at higher scales, which enhances spectral resolution. However, due to the uncertainty principle, this comes at the cost of worsened spatial localization compared to standard needlets. As a result, shrinking needlets are particularly suited for applications that prioritize fine frequency discrimination over spatial compactness, such as high-frequency asymptotic analysis.

By the technical point of way, shrinking needlets correspond to subexponential center sequences, that is 
\[
\underset{j \rightarrow \infty}{\lim} \frac{\log S_j}{j} = 0, 
\]
or equivalently,
\begin{equation*}
		\lim_{j\rightarrow \infty }h_{j}=  1,
\end{equation*}
see for additional details \cite{d25}.
In this case we can rewrite 
\[
h_j = 1+ \varepsilon_j,
\]
where the \emph{shift} $\varepsilon_j$ converges to zero as $j\rightarrow\infty$.
As shown in \cite{d25}, ensuring an appropriate balance between spectral and spatial localization requires the shift sequence to be sufficiently regular, with the divergence condition
\[\underset{j \rightarrow \infty}{\lim}\sum_{k=0}^{j-1} \varepsilon_j = \infty.
\]
To satisfy this, we consider shifts of the form
\begin{equation*}
\varepsilon_j = \frac{\gamma(j)}{j^p},
\end{equation*}
where $p \in \left(\left.0, 1\right]\right. $ and $\gamma: \mathbb{R}^+\to \mathbb{R}^+$ is a slowly varying function in the sense of Karamata (see \cite{bingham}), that is, for any $\tau >0$ 
\[
\underset{j\rightarrow \infty }{\lim} \frac{\gamma(\tau j)}{\gamma(j)}=1,
\] 
with some additional conditions for $\gamma$ at the critical value $p=1$ to be discussed below, see \cite{d25}). 
In this case we have
\begin{equation}\label{eqn:Sj}
	S_j =\begin{cases} \exp \left(\gamma(j)\frac{j^{1-p}}{1-p}\right) & p\in(0,1)\\
		\exp \left(\gamma(j)\log j\right) & p=1
		\end{cases}.
\end{equation}
For any $j \geq 1$, he discretized version of shrinking needlets is constructed using a set of cubature weights $\lbrace \lambda_{j,k}: k=1,\ldots,K_j  \rbrace$ and associated cubature points $\lbrace \xi_{j,k}:k=1,\ldots,K_j  \rbrace$ satisfying appropriate quadrature properties on the sphere. The resulting needlet functions are defined as
\begin{equation*}
	\psi_{j,k}(x) = \sqrt{\lambda_{j,k}} \sum_{\ell \in \left[S_{j-1},S_{j+1}\right]} b_j (\ell) Z_{\ell} \left(\langle x, \xi_{j,k}\rangle\right),
\end{equation*}
where $\langle x, \xi_{j,k}\rangle$ is the usual great circle distance on the sphere. Note that \[
\lambda_{j,k} \sim S_{j-1}^{-2}; \quad K_j  \sim S_{j-1}^{2},
\]
see \cite{dmt24,d25}.
Under these conditions, in \cite{d25} the following nearly-exponential localization property is established. For each $x,y \in \mathbb{S}^2$, let $d_{\mathbb{S}^2}(x,y)$ denote the great circle distance between the two spherical locations. For all $j \geq 1$, $x\in \mathbb{S}^{2}$ and
for all integers $M$, there exist two constants $C_{M}>0$ and $C^\ast_M>0$ such that%
	\begin{equation*}\label{eqn:loc}
	\begin{split}
		\left \vert	\psi_{j,k} (x) \right \vert 
		&	\leq\frac{4 C_M  	\Sigma_{j;p}}{\left(1+	\Sigma_{j;p}\Theta_{j,k}\right)^{M}}, ,\end{split}
\end{equation*}  
and 
\begin{equation}\label{eqn:locPhi}
	\left \vert  \Phi_j (x,y)\right \vert \leq \left(S_{j+1}^2-S_{j-1}^2\right)\frac{C^\ast_M  }{\left(1+\Sigma_{j,p}d_{\mathbb{S}^2}(x,y)\right)^M}.
\end{equation}
where
\begin{equation}\label{eqn:locfunct}
	\Sigma_{j;p} = \varepsilon_jS_j =\begin{cases} \gamma(j) /j^p \exp\left( \frac{j^{1-p}\gamma(j)}{1-p}\right) & p \in (0,1)\\ 
		\gamma(j) \exp\left( \log j \gamma(j)\right) & p =1, \quad \gamma(j) \geq O \left(\left( \log j\right)^{-1+\delta}\right)\\
		\eta \left(\log j \right)^{\eta-1}& p =1, \quad \gamma(j) = \frac{\eta}{\log j}, \quad \eta>1\\
		\exp\left( j \gamma(j) \right)& p =0, \quad \gamma(j): \underset{j \rightarrow \infty}{\lim} \gamma(j)  =0.
	\end{cases}.
	\end{equation}
and  $\Theta_{j,k}=d_{\mathbb{S}^2}\left(x,\xi_{j,k}\right) :=\arccos (\left\langle x,\xi_{j,k}\right\rangle )$. 
In the first case, the exponential factor behaves as an exponential of a sublinear power of $j$ and we van label this stretched-exponential growth. We refer to the second case as the power-low growth. Indeed, the asymptotic behavior of $\Sigma_{j;-1}$ splits into three regimes, depending on the behaviour of $\gamma$. Indeed, wheter $\gamma(j)=O\left(\left(\log j\right)^{-1+\delta}\right)$, $\gamma(j) \rightarrow c>0$ or $\gamma(j) \rightarrow \infty$, we have sub-polynomial, polynomial or super-polynomial growth respectively.
In the third case, we obtain a pure power of $\log j$ and it is labeled logarthmic growth, while the fourth corresponds to the critical region between shrinking and stable regimes.     
\\
As a direct consequence of the localization property, the following bound holds for the $L^p$-norms of the needlets (see \cite{d25}). 
	For $q\in \left[ 1,+\infty \right] $,
\begin{equation*}
\zeta_q \Sigma_{j;p}^{\left( 1%
	-\frac{2}{q}\right) } \leq 	\left\Vert \psi _{j,k}\right\Vert _{L^{q}\left( \mathbb{S}^{2}\right)
	}=\left( \int_{\mathbb{S}^{2}}{\left\vert \psi _{j,k}\left( x\right)
		\right\vert ^{q}dx}\right) ^{\frac{1}{q}}\leq \zeta_q \Sigma_{j;p}^{\left( 1%
		-\frac{2}{q}\right) } . \label{eq:pnorm}
\end{equation*}
\subsection{Introducing Spherical Poisson Needlets}
In \cite{bdmt24}, the spherical random wave $T_{\ell;t}(x)$ has been defined as a random spherical wave of frequency $\ell$, constructed by averaging the harmonic projector $Z_\ell(x, y)$ over random locations $y \in \mathbb{S}^2$ sampled from a Poisson process with intensity $\nu_t$:
\[
T_{\ell;t}(x)= \sqrt{\frac{4\pi}{\nu_t(2\ell+1)}}\int_{\mathbb{S}^2}Z_{\ell} (x,y)N_t\left(	\operatorname{d} y\right),\] 
where $\lbrace N_t(\cdot)\rbrace $ is a Poisson process on the sphere with governing intensity measure
\[ E\left[N_t(A)\right]=\nu_t\mathcal{L}(A) \quad \mbox{ for all } A\in\mathcal{B}\left( \mathbb{S}^2 \right).\]
Heuristically, the model in \cite{bdmt24} defines Poisson spherical random waves as random superpositions of spherical harmonics of degree $\ell$ centered at random points 
$y$ drawn from a Poisson point process on $\mathbb{S}^2$, with intensity $\nu_t$.
More rigorously,
\begin{definition}
	Let $\left(\Theta,\mathcal{A}, \rho \right)$ be a $\sigma$-finite measure space, with no amotic $\mu$. A collection of random variables $\lbrace N(A) : A \in \mathcal{A} \rbrace$ 
	taking values in $\mathbb{Z}^+ \cup \lbrace \infty \rbrace$, is a Poisson random measure $\Theta$ with intensity measure $\rho$ if the following two properties hold:
	\begin{enumerate}  
		\item $N(A)$ has Poisson distribution with control $\rho(A)$ for every $A \in \mathcal{A}$;
	\item  If $A_1,\ldots, A_n$ is a collection of pairwise disjoint elements in $\mathcal{A}$, then the measures $N(A_1), \ldots,N(A_n)$ are independent.
	\end{enumerate}
\end{definition}
In our setting, $\Theta = \mathbb{R}^+\times \mathbb{S}^2$, and $\mathcal{A}$ is the class of Borel subsets of $\Theta$, while $N$ is a Poisson random measure with homogeneous intensity is separable as \(\rho= \lambda \times \mathcal{L},\) multiplying a random measure taken as the time component, with the Lebesgue measure on the sphere. Regarding to the time component, assuming that $\lambda\left(\lbrace 0 \rbrace\right)=0$ and that $\lambda(\left[0,t\right])= \nu_t$ is monotonically increasing and diverging to infinity, we have that $t= \to \nu_t$ is the distribution function of $\lambda$. Thus, the mapping $A \to N_t(A)= [0,t]\times \mathcal{L} (A)$ is a Poisson random measure on $\mathbb{S}^2$, see \cite{dmp14,dmt24}, with the same distribution as 
\[
A \to \sum_{k=1}^{N_t(\mathbb{S}^2)} \delta_x(A),\]
where $\delta_x(\cdot)$ denotes the Dirac measure at the point $x$. Note that spherical Poisson waves can be rewritten as 
\[
T_{\ell;t}(x)= \frac{1}{\sqrt{\nu_t}} \sum_{i=1}^{N_t\left(\mathbb{S}^2\right)} Z_{\ell}\left(\langle x, \zeta_i\rangle\right),
\]
so that spherical Poisson random waves can be read as a weighted sum of a random number of deterministic waves, centred at random points $ \zeta_1,\zeta_2,\ldots$ uniformly distributed on the sphere. 
Using the addition formula we find a third expression for the waves,  
\[
T_{\ell;t} (x) = \sum_{m=-\ell}^{\ell} \hat{a}_{\ell,m}Y_{\ell,m}(x), 
\]
where the random harmonic coefficients are given by the formula
\[
 \hat{a}_{\ell,m} = \sqrt{\frac{4\pi}{(2\ell+1)\nu_t}} \sum_{i=1}^{N_t\left(\mathbb{S}^2\right)} Y_{\ell,m}\left(\zeta_{i} \right),
\]
see again \cite{dmt24}.\\
In this paper, we introduce an analogous construction that replaces spherical harmonics with properly normalized shrinking needlet projectors, motivated by their sharp spatial concentration, which enables localized and scale-sensitive analysis of spherical random fields.\\
More in detail, the spherical Poisson needlet with rate $\nu_t$  is defined by 
\begin{equation*}\label{eqn:PSW}
\Psi_{j;t}(x)= \frac{1}{\sqrt{\nu_t}\sigma_j}\int_{\mathbb{S}^2}\Phi_{j} (x,y )N_t\left(\operatorname{d} y\right),
\end{equation*} 
where the normalization factor $\sigma_j$ is defined by
\begin{equation*}\label{eqn:sigma}
\sigma_j^2 = \sum_{\ell \geq
	0}\frac{2\ell +1 }{4\pi} b^4_j\left( \ell\right),
\end{equation*}	
 such that 
 $$
 \mathbb{E} \left\{ \Psi_{j;t}(x)\right\} = 0; \quad \operatorname{Var} \left(\Psi_j(x)\right)= 1.$$
 Heuristically, the field $\Psi_{j;t}(x)$ is constructed by randomly superimposing localized shrinking needlet functions $\Phi_j(x, y)$, each centered at a random point $y$ on the sphere drawn from a Poisson process with intensity $\nu_t$; the normalization by $\sqrt{\nu_t}\sigma_j$ ensures that the resulting field has unit variance and captures local random fluctuations at scale $j$. The following proposition makes explicit the dependence of $\sigma_j^2$ on the scale parameter $j$, showing its asymptotic growth in terms of $S_j^2 \varepsilon_j$.  
 \begin{proposition}\label{prop:sigma}
 	Let the notation above prevail. Then it holds that
 	\[
 	\sigma^2_j \sim  \frac{C_b}{\pi} (S_j^2 \varepsilon_j),
 	\]
 	where $C_b$ is a positive real constant depending on the weight function $b$. 
 \end{proposition}
The proof of Proposition \ref{prop:sigma} is available in Section \ref{sec:proofs}.\\
An alternative notation, which more clearly highlights the intrinsic randomness in the sample size, is given by
\begin{equation*}
\Psi_{j;t}(x) = \frac{1}{\sqrt{\nu_t}\sigma_j}\sum_{i=1}^{N_t(\mathbb{S}^2)} \Phi_j (x,z_i),
\end{equation*}
as a sum of a random number of deterministic needlet waves, centred at points $z_i$ uniformly distributed on the sphere. \\
Finally, we can express $\Psi_{j;t}$ as a linear combination of random Poisson needlet coefficients weighted by the corresponding (deterministic) shrinking needlets, that is,
\begin{equation*}
	\Psi_{j;t}(x) = \sum_{k=1}^{K_j} \hat{\beta}_{j,k;t} \psi_{j,k} (x),
\end{equation*} 
where 
\[
\hat{\beta}_{j,k;t} = \frac{1}{\sqrt{\nu_t}\sigma_j} \int_{\mathbb{S}^2} \psi_{j,k}(x) N_t(\operatorname{d}x).
\]
Observe that 
\begin{equation*}
	\begin{split}
	&	\mathbb{E}\left[\hat{\beta}_{j,k;t} \right]=0;\\
	& \operatorname{Var}\left(\hat{\beta}_{j,k;t} \right)=\frac{1}{\sigma_j^2} \left \Vert \psi_{j,k}\right \Vert_2^2 = \frac{\bar{\sigma}^2}{S_j^2 \varepsilon_j},
	\end{split}
\end{equation*}
where 
\[
\bar{\sigma}^2=\frac{\pi}{C_b}  \left \Vert \psi_{j,k}\right \Vert_2^2,
\]
in view of Proposition \ref{prop:sigma}. Even if the definition can be misleading, $\bar{\sigma}$ does not depend either on $j$ or $k$. Indeed, in view of \eqref{eqn:locfunct}, $\left \Vert \psi_{j,k}\right \Vert_2$ is a positive real constant does not depend explicitly on $j,k$ but is a positive constant bounded by 1 (see also \cite{d25}).


\subsection{Key Estimates on Poisson–Wiener Chaoses}
We collect here several useful bounds for random variables that lie in the first Wiener–It\^o chaos of a Poisson measure. \
We will work with some probability metrics, each metrizing a topology strictly finer than weak convergence.  Throughout, for a function $g:\mathbb{R}^d\to\mathbb{R}$ we write
$$
\|g\|_{\mathrm{Lip}}
=\underset{x\in\mathbb{R}^d}{\sup}\|\nabla g(x)\|_{\mathbb{R}^d}.
$$
A straightforward generalization is given as follows. For any test function $g\in C^k(\mathbb{R}^d)$, we introduce its partial–derivative modulus
$$
M_k(g)=\sup_{x\neq y} \frac{\bigl\| D^{\,k-1}g(x)-D^{\,k-1}g(y)\bigr\|_{\mathrm{op}}}
{\|x-y\|},
$$
where $\|\cdot\|_{\mathrm{op}}$ is the operator norm on $(\mathbb{R}^d)^{\otimes(k-1)}$, defined by 
$$
\|T\|_{\mathrm{op}} := \sup \left\{ |T(v_1, \dots, v_{k-1})| : v_i \in \mathbb{R}^d,\|v_i\| =  1,i=1,\ldots,k-1 \right\},
$$
where $T$ is a multilinear form acting on $k-1$ vectors in $\mathbb{R}^d$.
In particular, for $g\in C^2(\mathbb{R}^d)$, we have that 
$$
M_2(g)
=\underset{x\in\mathbb{R}^d}{\sup}\|\mathrm{Hess}\,g(x)\|_{\mathrm{op}}.
$$
while, for $g\in C^3(\mathbb{R}^d)$, it holds that
\[
M_3(g) = \sup_{x \neq y} \frac{\left\| \mathrm{Hess}\, g(x) - \mathrm{Hess}\, g(y) \right\|_{\mathrm{op}}}{\|x - y\|}.
\]
Heuristically, $M_2(g)$ controls the curvature of the function $g$, i.e., how sharply $g$ can bend. This corresponds to a uniform bound on its Hessian, containing its second derivatives. In contrast, $M_3(g)$ measures how fast this curvature can change, that is, how much the Hessian varies across space. A bound on $M_3(g)$ ensures that the second derivatives of $g$ are Lipschitz continuous, and therefore that $g$ behaves smoothly at third order as well. 

\begin{definition}[Wasserstein and Kolmogorov distances]
	If $X,Y$ are $\mathbb{R}^d$-valued random vectors with $\mathbb{E}\|X\|<\infty$ and $\mathbb{E}\|Y\|<\infty$, then their Wasserstein distance is
	$$
	d_W(X,Y)
	= \underset{\|g\|_{\operatorname{Lip}\leq 1}}{\sup}
	\left \vert\mathbb{E}g(X)-\mathbb{E}g(Y)\right \vert.
	$$
	Also, if \(X, Y\) are real-valued random variables, their Kolmogorov distance is given by
	\[
	d_{\mathrm{Kol}}(X, Y)
	= \sup_{x \in \mathbb{R}} \left| \mathbb{P}(X \leq x) - \mathbb{P}(Y \leq x) \right|.
	\]
\end{definition}
\begin{definition}[$d_2$–distance]
	Under the same moment assumptions, define
	$$
	d_2(X,Y)
	=
	\underset{\substack{g\in C^2(\mathbb{R}^d)\\\|g\|_{\mathrm{Lip}}\le1,\;M_2(g)\le1}}{\sup}
	\bigl|\mathbb{E}\,g(X)\;-\;\mathbb{E}\,g(Y)\bigr|.
	$$
\end{definition}
\begin{definition}[$d_3$-distance] Under the same assumptions, 
	$$
	d_3(X, Y) \;=\;
	\sup_{\substack{g \in C^3(\mathbb{R}^d); \|g\|_{\mathrm{Lip}} \le 1 \\ M_2(g) \le 1; M_3(g) \le 1}}
	\bigl| \mathbb{E}[g(X)] - \mathbb{E}[g(Y)] \bigr|,
	$$
\end{definition}
The first two results (univariate and multivariate) are versions of the Fourth Moment Theorem for random variables belonging to the Wiener chaos of a general Poisson random measure proved in \cite{dvz} . The third is a functional convergence bound from \cite{bc}.
\subsubsection*{Poisson integrals and the first chaos}
Let $(\Omega,\mathcal{A},\mathbb{P})$ be a probability space supporting a Poisson random measure $N$ on a measurable space $(Z,\mathcal{Z})$ with intensity $\mu$.  Denote by
$$
\widehat N(\cdot)=N(\cdot)-\mu(\cdot)
$$
the compensated Poisson measure.  For $p\geq 1$, we write $L^p(Z)$ for $L^p(Z,\mathcal{Z},\mu)$.
\begin{definition}
	For any deterministic $h\in L^2(Z)$, define its first-order Wiener–It\^o integral
	$$
	I_1(h)=\int_Z h(z)\widehat N(\operatorname{d}z).
	$$
	The collection of all such $I_1(h)$, which is closed in $L^2(\Omega)$, is called the first Poisson Wiener chaos, denoted $\mathcal{C}_1$.
\end{definition}
\subsubsection*{Fourth–Moment bounds: univariate case}\label{sec:appuni}
Let $F\in\mathcal{C}_1$ satisfy $\mathbb{E}[F]=0$, $\mathrm{Var}(F)=1$, and $\mathbb{E}[F^4]<\infty$.  Write $Z\sim N(0,1)$.  Then there is a universal constant $\tilde{C}=\frac{1}{\sqrt{2\pi}}+\frac{2}{3}$ such that
\begin{equation}\label{eqn:wass}
d_W(F,Z)\leq \tilde{C}\sqrt{\mathbb{E}[F^4]-3},
\end{equation}
and moreover
\begin{equation}\label{eqn:kol}
d_{\mathrm{Kol}}(F,Z)\leq
\left( 11+ \sqrt{\mathbb{E}[F^4]}+(\mathbb{E}[F^4])^{1/4}
\right)
\sqrt{\mathbb{E}[F^4]-3},
\end{equation}
where $d_W$ and $d_{\mathrm{Kol}}$ are the Wasserstein and Kolmogorov distances respectively (cf. \cite{dvz}[Theorem 2.1 \& Corollary 1.3].
\subsubsection*{Fourth–Moment bounds: multidimensional case} Let $F_{(d)}=(F_1,\ldots,F_d)$, centered with covariance matrix $\mathcal{C}_d$, and $Z_d\sim \mathcal{N}\left(0,\mathcal{C}_d \right)$. Them it holds that 
\begin{equation}\label{thm:A3bis}
	\begin{split}	
		\left\vert \mathbb{E}[h(F_{(d)})] - \mathbb{E}[h(Z_d)]\right \vert
		& \leq B_3(h,d) \sum_{i=1}^{d}\sqrt{\mathbb{E}\left[F_i^4\right]-3\mathbb{E}^2\left[F_i^2\right]},
	\end{split}
\end{equation}
where
\begin{equation}\label{eqn:Btre}
B_3(h,d)=A_2(h;d) + \frac{2}{9}\sqrt{d \operatorname{Tr}\left(\mathcal{C}_d\right)}M_3(h), \quad A_2(h,d)=\frac{\sqrt{2d}}{4}M_2(h).
\end{equation}
cf. \cite{dvz}[19, Thm 1.7 \& Cor 1.8].\\
If in addition $\mathcal{C}_d$ is positive definite, then for any $h\in C^2(\mathbb{R}^d)$,
\begin{equation}\label{thm:A3}
	\begin{split}	
		\left \vert \mathbb{E}[h(F_{(d)})] - \mathbb{E}[h(Z_d)]\right \vert
		& \leq B_2(h,d) \sum_{i=1}^{d}\sqrt{\mathbb{E}\left[F_i^4\right]-3\mathbb{E}^2\left[F_i^2\right]},
	\end{split}
\end{equation}
where
$$
B_2(h,d)=A_1(h;d) + \frac{\sqrt{2\pi}}{6}\left \Vert \mathcal{C}^{-\frac{1}{2}}_d \right\Vert_{\operatorname{op}}\operatorname{Tr} \mathcal{C}_dM_2(h); \quad  A_1(h,d)=\frac{\left \Vert \mathcal{C}^{-\frac{1}{2}}_d \right\Vert_{\operatorname{op}}}{\sqrt{\pi}}M_1(h).
$$
cf. \cite{dvz}[Remark 4.3].

\subsubsection*{Functional convergence in Hilbert–Poisson Spaces}
Finally, one can lift these bounds to a functional setting (e.g. for processes indexed by a parameter set). A detailed statement and proof can be found in \cite{bc}; here, we only present the results as they have been adapted to our setting.\\
Let $K$ be a separable Hilbert space, and suppose $X$ is a $K$-valued random variable with finite second moment, $\mathbb{E}\|X\|_K^2<\infty$.  Its covariance operator $S:K\to K$ is defined by

$$
S u =\mathbb{E}\left[\langle X,u\rangle_K X\right],
$$
for every $u\in K$.  One checks that $S$ is self-adjoint, positive, and of trace-class, with
$$
\mathrm{Tr}S=\mathbb{E}\left[\|X\|_K^2\right].
$$
We work in the Banach space of all trace-class operators on $K$, endowed with the trace norm
$\|A\|_{\mathrm{Tr}}=\mathrm{Tr}\,\sqrt{A A^*}$.  A natural subspace is the Hilbert–Schmidt class $\mathrm{HS}(K)$, whose norm is
$\|A\|_{\mathrm{HS}}=\sqrt{\mathrm{Tr}(A A^*)}$.

Now assume in addition that $X$ lives in the first Poisson Wiener chaos, so it can be written as a single compensated Poisson integral, and has a finite fourth moment, $\mathbb{E}\|X\|_K^4<\infty$.  Let $Z$ be a centered Gaussian random element in $K$ whose covariance operator is also $S$.  We compare the laws of $X$ and $Z$ via a smooth-distance metric $d_3$, which controls differences of expectations against three-times differentiable test functions on $K$.
Under the above hypotheses,
$$
d_3\left(X,Z\right)\lesssim
\sqrt{
	\left(\mathbb{E}\|X\|_K^4 - \mathbb{E}\|X\|_K^2\right)^2-
	2\|S\|_{\mathrm{HS}}}.
$$
More precisely, there exists an absolute constant $C$ so that
$$
d_3(X,Z)\leq
\left(\frac{1}{4}+\frac{1}{2}\sqrt{\mathbb{E}\|X\|_K^2} \right)\sqrt{
	\left(\mathbb{E}\|X\|_K^4 - \mathbb{E}\|X\|_K^2\right)^2-
	2\|S\|_{\mathrm{HS}}}.
$$
see \cite{bc}. This bound quantifies how close the Poisson-chaos field $X$ is to its Gaussian counterpart $Z$ in a functional form, by measuring both the fourth-moment deficit and the Hilbert–Schmidt size of the covariance.

\section{Main results}\label{sec:main}
In this section, we present our main quantitative central limit theorems in three complementary settings.  We begin by analyzing the Poisson needlet coefficients, establishing precise Gaussian approximations and convergence rates; for a detailed comparison with the corresponding results in the standard needlet framework, see \cite{dmp14}.  Next, we turn to finite‑dimensional distributions, proving a multivariate central limit theorem for the vectors $\{\Psi_j(x_i)\}_{i=1}^d$ of Poisson needlets evaluated at a fixed collection of well‑separated points on $\mathbb{S}^2$.  Finally, we address functional convergence, first in the $L^2$ topology and then in stronger Sobolev spaces, showing how these results complement the harmonic analysis analogues of \cite{bdmt24}.

\subsection{Central limit theorems for  Poisson needlet coefficients}\label{sec:coeff}
In this section, we develop quantitative central limit theorems for Poisson-driven needlet coefficients on the sphere, leveraging the powerful Peccati–Zheng bounds in the Poisson Wiener chaos \cite{peczen}, as in \cite{dmp14} for the standard needlet cases. Indeed, by combining these Malliavin–Stein techniques with the sharp localization of needlets, we derive explicit rates for the Gaussian approximation in both smooth‐test and Wasserstein distances, expressed in terms of the needlet resolution and  the time component of Poisson intensity.\\
In what follows, in the main propositions, we will establish bounds for the coefficients in a format that is particularly well-suited for comparison with finite-dimensional distributions, both in the unidimensional and multidimensional settings. In the corresponding corollaries, the focus shifts toward elucidating the convergence criteria for these coefficients, providing a clearer understanding of the conditions under which convergence holds.
\subsubsection*{Bounds in dimension one}
First of all, 
  $Z \sim \mathcal{N}(0,\sigma^2_Z)$ indicates a centered Gaussian random variable with variance  $\frac{\bar{\sigma}^2}{S_j^2 \varepsilon_j}$. 
\begin{proposition}\label{prop:dimone0}
	Let $\tilde{C}=\frac{1}{\sqrt{2\pi}}+\frac{2}{3}$. For every $j \geq 1$,$k=1,\ldots, K_j$ and every $t>0$, it holds that
	\[d_{\operatorname{W}} \left(\hat{\beta}_{j,k;t},Z \right) \leq \tilde{C} {\bar{\sigma}^4\zeta_4^4}  \frac{1}{S_j\sqrt{\nu_t}}. 
	\]
	\end{proposition}
	The proof is available in Section \ref{sec:proofs}.\\
	For any $j \geq 1$, and $k=1,\ldots,K_j$, we define the \emph{renormalized Poisson needlet coefficient}:
		\begin{equation}\label{eqn:normal}
				\tilde{\psi}_{j,k;t} = \frac{1}{\sqrt{\nu_t}\tilde{\sigma}^2}\int_{\mathbb{S}^2}\psi_{j,k}(x)N_t(d x),
		\end{equation}
		where $\tilde\sigma ^2 =\left\Vert\psi_{j,k} \right\Vert_2^2.$
Note that the renormalized coefficient
$\tilde{\psi}_{j,k;t}$ can be written in terms of the original (global unit-variance) coefficient $\hat{\psi}_{j,k;t}$ via the scaling factor
$$
\frac{S_j\,\sqrt{\varepsilon_j}}{\bar\sigma}.
$$
Since $S_j\to\infty$ and $\varepsilon_j\to0$ as $j\to\infty$, this factor captures precisely how the shrinking‐needlet normalization (through $S_j$) and the small‐scale variance (through $\sqrt{\varepsilon_j}$) combine to yield the constant limit variance $\bar\sigma^2$. In other words, as $j$ grows, the product $S_j\sqrt{\varepsilon_j}$ compensates exactly for the shrinking support of the needlets so that $\tilde\psi_{j,k;t}$ remains on a common scale.	
\begin{corollary}\label{cor:dimone0}
	Let the notation above prevail. For every $j \geq 1$, $k=1,\ldots, K_j$ and every $t>0$, it holds that
	\[d_{\operatorname{W}} \left(\tilde{\beta}_{j,k;t},\tilde{Z} \right) \leq \tilde{C} \frac{\zeta_4^4}{\tilde{\sigma}^4} \sqrt{ \frac{S_j^2\varepsilon_j^2}{\nu_t}}. 
	\]
	where $\tilde{Z}\sim\mathcal{N}(0,1)$.
	It follows that, if $S_j=o(\nu_t^{1/2})$, $\tilde B_{j;t}$ converges in law to $\tilde Z$ 	
\end{corollary}
The proof is included in Section \ref{sec:proofs}.
\begin{remark}[Bound for Kolmogorov distance]
	By using Equation \eqref{eqn:kol}, we can show that
	\begin{equation*}
	d_{\mathrm{Kol}}\left(\tilde{\beta}_{j,k;t},\tilde{Z} \right)\leq
	C \left( 11+ \sqrt{3+   \frac{S_j^2\varepsilon_j^2}{\nu_t}}+o\left(\frac{S_j^2\varepsilon_j^2}{\nu_t}
	\right)\right)
 \sqrt{ \frac{S_j^2\varepsilon_j^2}{\nu_t}} = O\left( \sqrt{ \frac{S_j^2\varepsilon_j^2}{\nu_t}}\right).
	\end{equation*} 
\end{remark}
\begin{remark}[Convergence via Scaling Structure]\label{rem:conv1}
	In this remark, we link the convergence rate of our quantitative CLT to the specific scaling behavior of $S_j$ , and we derive an explicit condition on the relationship between the resolution index $j$ and the Poisson intensity $\nu_t$,
	ensuring convergence in law.\\
	\begin{itemize}
		\item \emph{Stretched exponential case}.
	In the stretched exponential case, we can choose 
	\[
	j(t) = o\left(\left(\log \nu_t\right)^{\frac{1}{1-p}}\right).
	\]
 Indeed, since $\gamma$ is slowly varying, it follows that 
 \[
 j^{1-p}(t) \gamma(j(t)) = o\left(\log \nu_t\right),
 \]
and hence 
\[
\left(\Sigma_{j(t);p} \right)^2 = \gamma(j(t))^2 j^{-2p}(t)\exp\left(\frac{2 j^{1-p}(t)\gamma(j(t))}{1-p}\right) = o (\nu_t).
\]
	\item \emph{Power law case}. In the power-law case, 
$$
\Sigma_{j;-1}  =\gamma(j)\exp\left(\log j\gamma(j)\right) =\gamma(j) j^{\gamma(j)},
$$
so that
$$
\left(\Sigma_{j;-1}\right)^2
=\gamma(j)^2 j^{2\gamma(j)}.
$$
Requiring $\gamma(j)^2\,j^{2\gamma(j)} = o\left(\nu_t\right)$ is equivalent to
$$
2\log\gamma(j) + 2 \gamma(j) \log j
= o\left(\log\nu_t\right).
$$
The first addend $2\log\gamma(j)$ is at most logarithmic in $j$, thus the dominant condition is
\[
\gamma(j) \log j = o\left(\log\nu_t\right).
\]
Considering for example the pure polynomial growth, that is $\gamma \rightarrow c>0$, 
this becomes $c \log j(t)=o(\log\nu_t)$, i.e.\ $j(t)=\nu_t^{o(1)}$.\\
	\item \emph{Logarithmic case}.
For the logarithmic growth, the convergence condition reads
$$
\eta(\log j)^{\eta-1}
= o\left(\nu_t\right),
$$
equivalent to
$$
\left( \eta -1 \right) \log \log j)
= o\left( \log \nu_t\right),
$$
as for example 
$$
\log j(t)
=  \exp\left(o(\nu_t^{1/(\eta-1)})\right).
$$
Finally, for the critical growth, since $\gamma \rightarrow 0$ it suffices to choose 
$$ j(t)=o\left(\log \nu_t\right).$$
\end{itemize}
Table \ref{tab:j_vs_nu} collects all the conditions described in this remark.
\begin{table}[h]
	\centering
	\begin{tabular}{ll}
		\hline
		\vspace{1pt}
		\textbf{Case} &  \textbf{Growth of \(\displaystyle j(t)\)} \\ 
		\hline
			\vspace{1pt}	
		Stretched exponential \((0<p<1)\) & \(\displaystyle j(t)=o\left(\left(\log \nu_t\right)^{\frac{1}{1-p}}\right)\) \\[1ex]
		Power law \((p=-1)\) & \(\displaystyle j(t)=\nu_t^{o(1)}\) \\[1ex]
		Logarithmic growth \((\gamma(j)=\eta\log j,\;\eta>1)\) & \(\displaystyle \log j(t)=\exp\left(o\left(\nu_t^{\tfrac{1}{\eta-1}}\right)\right)\) \\[1ex]
		Critical case \((\gamma(j)\to0)\) & \(\displaystyle j(t)=o(\log \nu_t)\) \\
		\hline
	\end{tabular}
	\vspace{3pt}
	\caption{Growth‐rate conditions on \(j(t)\) ensuring \(\Sigma_{j(t);p}^2=o(\nu_t)\) in the various regimes for needlet coefficients in the unidimensional case.}
	\label{tab:j_vs_nu}
\end{table}
\end{remark}
\subsubsection*{Multidimensional bounds}
For any fixed $j\geq 1$, let us define preliminarily the following random vector of  dimension $K_j$ 
\begin{equation*}
	B_{j;t} =\left(\beta_{j,1;t}, \ldots,\beta_{j,K_j;t}\right)^T. 
\end{equation*}
Let $\Upsilon_{j;t}$ be the related $K_j \times K_j$ covariance matrix. Observe that each entry of this matrix corresponds to
	\begin{equation*}
	\Upsilon_{j;t} \left(k_1,k_2\right) = \mathbb{E} \left(\beta_{j,k_1;t} \beta_{j,k_2;t}\right) 
	\leq \frac{\bar{\sigma}^2 C_M}{S_j^2\varepsilon_j \left(1+ \Sigma_{j;p} d_{\mathbb{S}^2} (\xi_{j,k_1}\xi_{j,k_2})\right)^M }, 
\end{equation*}
for $\quad k_1,k_2=1,\ldots K_j$, see Lemma 5.1 in \cite{dmp14}. Choosing sufficiently sparse cubature points, as for instance such that \[
d_{\mathbb{S}^2} \left(\xi_{j,k_1}\xi_{j,k_2}\right) \sim  \Sigma_{j,p}^{-\frac{1}{2}},\] 
we obtain a covariance matrix converging to a diagonal one, with elements on the main diagonal given by
	\begin{equation*}
	\Upsilon_{j;t} \left(k,k\right) = \frac{\bar{\sigma}^2 }{S_j^2\varepsilon_j}.
\end{equation*}
Let now $Z_{K_j}$ be  centered Gaussian vector with covariance matrix $\Upsilon_{j;t}$.
\begin{proposition}\label{prop:multicoef}
 Let the notation above prevail. Then it holds that
 \begin{equation*}
 	d_3\left(B_{j;t},Z_{K_j}\right) \leq \underset{g \in \mathcal{C}^3}{\sup} B_3(g,S_j) \sqrt{\frac{S_j^2}{\nu_t}},
 \end{equation*}
 where 
 \[
 B_3(g,S_j) = \frac{\sqrt{2}S_j} {4}M_2(g) + \frac{2\bar \sigma ^2}{9} \frac{S_j}{\sqrt{\epsilon_j }}.
 \]
Note that \[d_3\left(B_{j;t},Z_{K_j}\right) = O\left(\frac{S_j^2}{\nu_t^{\frac{1}{2}}}\right).\]
Also, if $S_j = o\left(\nu_t^{\frac{1}{4}}\right)$, then $B_{j;t}$ converges in law to $Z_{K_j}$.
\end{proposition}
The proof is included in Section \ref{sec:proofs}.\\
Let now $\tilde B_{j;t}$, be the normalized random vector whose components are given by Equation \eqref{eqn:normal}
\[\tilde{\beta}_{j,k;t}=\frac{1}{\tilde{\sigma}\sqrt{\nu_t}}\int_{\mathbb{S}^2} \psi_{j,k}(x) N_t{\operatorname{d}x}; \quad k=1,\ldots K_j,\]
which is centered and with covariance matrix $\tilde \Upsilon$ of dimension $K_j \times K_j$. The entries of this matrix out of the main diagonal are bounded by
\[
\begin{split}
\left\vert 	\mathbb{E} \left[\tilde{\beta}_{j,k_1;t}\tilde{\beta}_{j,k_2;t}\right] \right \vert &= \frac{1}{\tilde{\sigma}^2} \int_{\mathbb{S}^2}\left \vert  \psi_{j,k_1}(x)\right \vert \left \vert\psi_{j,k_2}(x)\right \vert \operatorname{d}x \\
&= \frac{C_M\zeta_2}{\tilde{\sigma}^2 \left(1+ \Sigma_{j;p} d_{\mathbb{S}^2}\left(\xi_{j,k_1},\xi_{j,k_2}\right) \right)^M } ,
\end{split}
\]
in view of Lemma 5.1 in \cite{dmp14} properly adapted to flexible bandwidth needlets. Then, using a properly sparse choice of cubatore points, as $j \rightarrow \infty$, $\tilde \Upsilon$ converges to the identity matrix of dimension $K_j$, so that the following result holds. 
	\begin{corollary}\label{cor:multicoef}
		Let the notation above prevail, and 
Let $\tilde  Z_{K_j} \sim \mathcal{N}(0,1)$  Then it holds that
		\begin{equation*}
			d_3\left(\tilde B_{j;t},\tilde Z_{K_j}\right) \leq \underset{g \in \mathcal{C}^3}{\sup}\tilde B_3(g,S_j)\sqrt{ \frac{S_j^6\varepsilon_j^2}{\nu_t}},
		\end{equation*}
	$$
\tilde B_3(g,S_j)=\frac{\sqrt{2}S_j}{4}M_2(g) + \frac{2}{9}S_j^2M_3(g).$$
Since 
\[
d_3\left(\tilde B_{j;t},\tilde Z_{K_j}\right) \leq \left( \frac{S_j^5}{\nu_t^{\frac{1}{2}}}\right)
\]
Also, $\tilde B_{j;t}$ converges in law to $\tilde Z_{K_j}$ if $S_j=\nu_t^\frac{1}{10}$.	\end{corollary}
The proof is available in Section \ref{sec:proofs}.
\begin{remark}[Bounds for the distance $d_2$]\label{rem:dtwo}
		Note that, since $\Upsilon_{j;t}$ is positive definite, and noting that 
		$$
		\left\Vert \tilde \Upsilon_{j;t} ^{-\frac{1}{2}}\right\Vert_{\operatorname{op}} \simeq \sigma_j \simeq S_j\sqrt{\varepsilon_j} ; \quad \operatorname{Tr} \tilde \Upsilon_{j;t}  = \frac{K_j \bar{\sigma}^2}{S_j^2 \epsilon_j} =  \frac{ \bar{\sigma}^2}{\epsilon_j}
		$$
		we can get the following bound:
		\begin{equation*}
			d_2\left(B_{j;t},Z_{K_j}\right) \leq \underset{g \in \mathcal{C}^2}{\sup} B_2(g,S_j) \sqrt{\frac{S_j^2}{\nu_t}},
		\end{equation*}
		where 
		$$
		B_2(g,S_j)=\frac{S_j \sqrt{\varepsilon_j}}{\sqrt{\pi}}M_1(g)+ \frac{\sqrt{2\pi} S_j}{6\sqrt{\varepsilon_j}}M_2(g).
		$$
		It follows that 
		\[
			d_2\left(B_{j;t},Z_{K_j}\right) \leq O\left(\frac{S^2_j}{\sqrt{\varepsilon_j\nu_t}}\right),
		\]
		so that if $S_j = o\left(\nu_t^{\frac{1}{4}}\right)$, $B_{j;t}$ converges in law to $Z_{K_j}$.\\ 
		In the normalized case, we obtain
		\begin{equation*}\label{eqn:d2bound}
			d_2\left(\tilde B_{j;t},\tilde Z_{K_j}\right) \leq \underset{g \in \mathcal{C}^2}{\sup} \tilde B_2(g,S_j) \sqrt{\frac{S_j^6 \varepsilon_j^2}{\nu_t}},
		\end{equation*}
		where
		$$
		\tilde B_2(g,S_j)=\frac{1}{\sqrt{\pi}}M_1(g)+ \frac{\sqrt{2\pi}S_j^2}{6}M_2(g).
		$$
		It follows that, if $S_j=o(\nu_t^{1/10})$, $\tilde B_{j;t}$ converges in law to $\tilde Z_{K_j}$ 
		\end{remark}	
\begin{remark}[Comparison with existing quantitative central limit theorems for spherical Poisson needlets]
In \cite{dmp14} a set of quantitative central limit theorems have been established for Poisson coefficients on the sphere by means of quantitative central limit theorems developed by \cite{peczen}. More in detail, for any deterministic and regular kernels yielding to normalized Poisson integral  $F$,  the Wasserstein distance between the law of $F$ and the standard normal random variable $Z\sim \mathcal{N}(0,1)$.satisfies	
	$$
	d_W\left(F,Z\right)
	\leq
	\left \vert1 -\|h\|_{L^2(\mu_t)}\right\vert
	+
	\frac{1}{\|h\|_{L^2(\mu_t)}^3}
	\int_{S^2}\left \vert h(z)\right \vert^3\;\mu_t(\operatorname{d}z).
	$$
The results obtained here are consistent with those previously established, once the notation is appropriately adapted.	Indeed, for every $j \geq 1$,$k=1,\ldots, K_j$ and every $t>0$, in the normalized case, $h=\frac{1}{\sqrt{\varepsilon_j}\bar{\sigma}}\psi_{j,k}$; as in \cite{dmp14}, we obtain that
	\[
	\begin{split}
		d_{\operatorname{W}} \left( \beta_{j,k;t},Z \right)  & \leq  
= 	\frac{1}{\sqrt{\nu_t}\bar\sigma^3} \left \Vert \psi_{j,k}\right \Vert_3^3 \\
	&=O\left(\sqrt{\frac{S_j^2\varepsilon_j^2}{\nu_t} }\right) ,
	\end{split}
	\]
	as in Corollary \ref{cor:dimone0}. \\
In the multidimensional case, for fixed \(D\geq1\), and \(Y\sim N_D(0,C)\) be a centered Gaussian vector in \(\mathbb{R}^D\) with positive–definite covariance matrix \(C\).  For each \(t>0\), consider the random vector
\[
F_t=(F_{t,1},\dots,F_{t,D})^T
=
\left(\widehat N_t(h_{t,1}),\dots,\widehat N_t(h_{t,D})\right)^T,
\]
where each \(h_{t,a}\in L^2(S^2,\mu_t)\).  Denote by \(\Gamma_t\) the covariance matrix of \(F_t\), i.e.
\[
\Gamma_t(a,b)
=
\mathbb{E}\left[\widehat N_t(h_{t,a})\,\widehat N_t(h_{t,b})\right]
\;=\;
\langle h_{t,a},\,h_{t,b}\rangle_{L^2(\mathbb{S}^2,\mu_t)},
\quad a,b=1,\dots,D.
\]
Then the \(d_2\)–distance between the laws of \(F_t\) and \(Y\) satisfies
\begin{equation*}
	\begin{aligned}
		d_2(F_t,Y)
		&\leq
		\bigl\|C^{-1}\bigr\|_{\mathrm{op}}\,
		\bigl\|C\bigr\|_{\mathrm{op}}^{1/2}\,
		\bigl\|C-\Gamma_t\bigr\|_{\mathrm{HS}}
		\\&\quad
		{}+\;\sqrt{\tfrac{2\pi}{8}}\,
		\bigl\|C^{-1}\bigr\|_{\mathrm{op}}^{3}\,
		\bigl\|C\bigr\|_{\mathrm{op}}\,
		\sum_{i,j,k=1}^D
		\int_{\mathbb{S}^2}
		\bigl|h_{t,i}(x)\bigr|\,
		\bigl|h_{t,j}(x)\bigr|\,
		\bigl|h_{t,k}(x)\bigr|\,
		\mu_t(\operatorname{d}x)\\
		&=\Delta_1+\Delta_2.
	\end{aligned}
\end{equation*}
Moreover, by applying H\"older’s inequality one also gets the simplified bound
\begin{equation*}
	\begin{aligned}
		d_2(F_t,Y)
		&\leq
		\bigl\|C^{-1}\bigr\|_{\mathrm{op}}\,
		\bigl\|C\bigr\|_{\mathrm{op}}^{1/2}\,
		\bigl\|C-\Gamma_t\bigr\|_{\mathrm{HS}}
		\\&\quad
		+D^2\sqrt{\tfrac{2\pi}{8}}\,
		\bigl\|C^{-1}\bigr\|_{\mathrm{op}}^{3}\,
		\bigl\|C\bigr\|_{\mathrm{op}}\,
		\sum_{i=1}^D
		\int_{\mathbb{S}^2}
		\bigl|h_{t,i}(x)\bigr|^3\,
		\mu_t(\operatorname{d}x)\\
		&=\Delta_1+\Delta_2^\prime.
	\end{aligned}
\end{equation*}
To set the comparison with \cite{dmp14}, here we consider $\tilde B^\prime_{j;d,t}$, vectors of needlet coefficients of dimension $d$ first fixed and then set up to grow, aimed to be simultaneously evaluated. In both the cases, to control $\Delta_1$, analogously to \cite{dmp14}, we assume that the cubature points at which these coeﬃcients are evaluated satisfy the condition:
\[
\underset{k_1 \neq k_2=1,\ldots,d_t}{\inf}d_{\mathbb{S}^2}\left(\xi_{j,k_1},\xi_{j,k_2} \right)\sim \Sigma_{j,p}^{-\frac{1}{2}}.
\]
In this case, $\Delta_1$ can be made as small as we aim just by choosing properly $M$.
Consider now $\Delta_1$, with $D=K_j=S_j^2$. It is straightforward to see that 
\[
\Delta_1 \leq O\left( \frac{S_j^7 \varepsilon_j}{\sqrt{\nu_t}}\right),
\]
which is slower that the bound established in Corollary \ref{cor:multicoef}. It could anyway be possible to improve this result, as in \cite{dmp14}[Theorem 5.5]), where the authors have used the concentration properties to show that $\Delta^{\prime } = o(\Delta)$.

\end{remark}
\begin{remark}[Convergence via Scaling Structure in the multidimensional case]
The regularity conditions on $j$ carry over unchanged to the multidimensional setting, since in both cases they arise from balancing $\log S_j$ against a constant multiple of $\log\nu_t$ (namely $\tfrac12\log\nu_t$ versus $\tfrac1{10}\log\nu_t$), and hence, up to that constant factor, the required relationship between $j$ and $\nu_t$ remains the same stated in Remark \ref{rem:conv1} and collected in Table \ref{tab:j_vs_nu}.
\end{remark}
\subsection{Finite dimensional distributions}\label{sec:fdd}
In this section, we study quantitative central limit theorems for finite-dimensional vectors formed by spherical Poisson needlets $\lbrace \Psi_j(x_i): i = 1,\ldots,d\rbrace$ computed over a fixed set of $d$ spherical locations that are sufficiently distant from each other. These random vectors arise from projecting a Poisson random measure on the sphere onto a localized needlet frame at resolution level $j$. Our aim here is to establish explicit bounds on the distance between the law of suitably normalized $\lbrace \Psi_j(x_i): i = 1,\ldots,d\rbrace$ and a multivariate Gaussian distribution, with a proper choice of metrics. 
\subsubsection*{Unidimensional case}
We begin by analyzing the one-dimensional case, both for simplicity of exposition and because the arguments naturally extend to arbitrary finite dimension $d$.
\begin{theorem}\label{thm:fin1d}
Let the notation above prevail. Then it holds that
\[
d_W\left(\Psi_{j;t}(x),Z\right)\leq \tilde{c} \sqrt{\frac{S_{j}^4 \varepsilon^2_j}{\nu_t}} ,
\]
where $\tilde{c} = \left(\frac{1}{\sqrt{2\pi}}+\frac{2}{3}\right)C_\infty$. Also, $\Psi_{j;t}(x)$ converges in law to $Z$ if $S_j=o\left(\nu_t^{\frac{1}{4}}\right)$.
\end{theorem}
The proof is available in Section \ref{sec:proofs} 
\begin{remark}[Bound for Kolmogorov distance for finite dimensional distributions]
Also in this case Equation \eqref{eqn:kol} can be used to derive a quantitative central limit theorem in the Kolmogorov distance rather than the Wasserstein distance: 
	\begin{equation*}
		\begin{split}
		d_{\mathrm{Kol}}\left(\Psi_{j;t}(x),Z\right) & \leq
		C_\infty \left( 11+ \sqrt{3+   \frac{S_j^4\varepsilon_j^2}{\nu_t}}+o\left(\frac{S_j^4\varepsilon_j^2}{\nu_t}
		\right)\right)
		\sqrt{ \frac{S_j^4\varepsilon_j^2}{\nu_t}} \\ 
		&= O\left( \sqrt{ \frac{S_j^4\varepsilon_j^2}{\nu_t}}\right),
		\end{split}
	\end{equation*} 
	 if $S_j=o\left(\nu_t^{\frac{1}{4}}\right)$.
\end{remark}
\subsubsection*{Multidimensional case}
We now turn to the multidimensional setting, extending the results established in the one-dimensional case to vectors of Poisson needlet evaluated over a fixed collection of $d$ well separated locations on the sphere. While the technical structure becomes richer in higher dimensions, the core probabilistic and analytical ingredients remain essentially unchanged. In particular, the arguments used in the univariate setting generalize naturally, allowing us to derive quantitative multivariate central limit theorems with explicit control on the dependence in $d$.\\
More in detail, to handle joint distributions at finitely many locations, fix distinct points $x_1,\dots,x_d\in \mathbb{S}^2$ and set
$$
\Xi_{d;j,t}=\left(\Psi_{j,t}(x_1),\dots,\Psi_{j,t}(x_d)\right).
$$
For $i_1,i_2=1,\dots,d$, its covariance matrix $\Gamma_{d;j,t}$ has entries
\[
\begin{split}
\Gamma_{d;j,t}\left(i_1,i_2\right)
=\mathbb{E}\left[\Psi_{j,t}(x_{i_1})\,\Psi_{j,t}(x_{i_2})\right]
=\frac{1}{\sigma_j^2}\sum_{\ell = S_{j-1}}^{S_{j+1}} b^{4}(\ell)Z_\ell\left(\langle x_{i_1},x_{i_2}\rangle\right),
\end{split}
\]
so in particular for each diagonal entry we have that $\Gamma_{d;j,t}\left(i,i\right)=1$, $i=1,\ldots,d$.
To ensure that the needlet covariance matrix $\Gamma_d$ converges to the identity as $j \to \infty$, it is sufficient to impose a minimal separation between the sample points on the sphere. Specifically, assume that the points $\{x_1, \dots, x_d\} \subset \mathbb{S}^2$ satisfy
\begin{equation*}
	\min_{1 \leq i < k \leq d} \, d_{\mathbb{S}^2}(x_i, x_k) \;\geq\; \delta ,
\end{equation*}
where $d_{\mathbb{S}^2}(\cdot, \cdot)$ denotes the geodesic distance on the sphere. Under this condition, standard localization estimates for needlets imply that for every $M > 0$ there exists a constant $C_M$ such that for $i \ne k$,
\begin{equation*}
	\left| \Gamma_{d;\, i, k} \right| \;\leq\; \frac{C_N}{(1 + \Sigma_{j;p} \delta)^N},
\end{equation*}
which tends to zero faster than any polynomial rate as $j \to \infty$, since $S_j \to \infty$. Hence, $\Gamma_d \to I_d$ entrywise as $j \to \infty$. Also, $\operatorname{Tr}(\Gamma_d) \rightarrow d$.
\\
Before stating the next theorem, we provide the reader with a quick note on notation, following \cite{bdmt24}.  In our definition of $d_3(\cdot,\cdot)$, one could restrict the supremum to test functions $g\in C^{2}$ whose derivatives satisfy $M_2(g)\le1$ and $M_3(g)\le1$, in analogy with the usual Wasserstein distance (where one takes Lipschitz functions of constant at most one).  Recall that, for $k\in\mathbb N$ and $g\in C^{k-1}(\mathbb R^d)$, we set
$$
M_k(g)
=\sup_{x\neq y}\frac{\|D^{\,k-1}g(x)-D^{\,k-1}g(y)\|_{\mathrm{op}}}{\|x-y\|_2},
$$
where $\|\cdot\|_{\mathrm{op}}$ is the operator norm on multilinear forms.  We also write $\|x\|_2$ for the Euclidean norm of $x\in\mathbb R^d$ and $\|A\|_{\mathrm{op}}$ for the induced operator norm of a matrix $A$.  Finally, $Z_d$ denotes the centered Gaussian vector in $\mathbb R^d$ with covariance matrix $\Gamma_d$.
\begin{theorem}\label{thm:multi}
Let $\Xi_{d;j,t} = \left(\Psi_{j}(x_1),\dots,\Psi_j (x_d)\right)$
be the $d$-dimensional vector of needlet spherical Poisson waves, all belonging to the first‐order Wiener-Poisson chaos, and let $Z_d\sim N_d(0,\Gamma_d)$ be the centered Gaussian vector with the same covariance $\Gamma_d$.  Then, as $j\to\infty$,
\[
\begin{split}
d_3\left(\Xi_{d;j,t} ,Z_d\right)\leq
& 
\sup_{g\in C^3}B_3(g;d) d \sqrt{\frac{S_{j}^4 \varepsilon^2_j}{\nu_t}}
\end{split}
\]
where
\begin{equation*}
	B_3(g,d)=\frac{\sqrt{2d}}{4}M_2(g)+ \frac{2d}{9}M_3(g).
\end{equation*}
Also, $\Xi_{d;j,t}$ converges in law to $Z_d$ if $S_j=o\left(\nu_t^{\frac{1}{4}}\right)$.
\end{theorem}
The proof is deferred to Section \ref{sec:proofs}.
\begin{remark}[Convergence via Scaling Structure for finite dimensional distributions]
Also in this case, the regularity conditions on $j$ remain unchanged for finite‑dimensional vectors, since in every case they follow from balancing $\log S_j$ against a fixed multiple of $\log\nu_t$ (e.g.\ $\frac{1}{2}\log\nu_t$ versus $\frac{1}{4}\log\nu_t$); consequently, up to that constant factor, the relationship between $j$ and $\nu_t$ is exactly as stated in Remark \ref{rem:conv1} and summarizad in Table \ref{tab:j_vs_nu}.
\end{remark}
\begin{remark}[Linking Finite-Dimensional Convergence and the Random Needlet Coefficients]
	As in \cite{bdmt24}, we want to link the results between the quantitative bounds for finite dimensional distributions and the results concerning the needlet coefficients. First, observe that using the typical quadrature formula for needlets (see \cite{npw1,dmt24}), we can discretize each component of $\Xi_{j;t}$ as
		\[
	\Psi_{j,t}(x_i)\;=\;\sum_{k=1}^{K_j}\hat \beta_{j,k;t}\psi_{j,k}(x_i),
	\qquad i=1,\dots,d.
	\]
	We can then introduce the deterministic $d \times K_j$ evaluation matrix
	$$
	R_{j,d}
	=\left[\psi_{j,k}(x_i)\right]_{1\leq i\le d,1\leq k\leq K_j}, 
	$$
	which is bounded by $ \tilde C_M dS_j\varepsilon_j$, by means of the localization property \eqref{eqn:locfunct},
	and stack the needlet‐coefficients into the random vector
	$$
	B_{j;t}
	=\left(\hat\beta_{j,1;t},\dots,\hat\beta_{j,K_j;t}\right)^T.
	$$	
	Then it holds that 
	$$
	\Xi_{d;j,t}  =R_{j,d} B_{j;t}.
	$$	
	Since $R_{j,d}:\mathbb{R}^{K_j}\mapsto\mathbb{R}^d$ is a fixed linear map, any test‐function
	$g:\mathbb{R}^d\mapsto\mathbb{R}$ lifts to
	$\tilde g:=g\circ R_{j,d}\colon \mathbb{R}^{K_j}\mapsto\mathbb{R},$
	and the fourth‐moment bound on $B_{j;t}$ yields
	\[
	\begin{split}
	\left \vert \mathbb{E} \left[g(\Xi_{;j,t})\right]-\mathbb{E}\left[g(Z_d)\right]\right \vert
	&=	\left \vert \mathbb{E} \left[\tilde g(B_{j;t})	\right]-\mathbb{E}\left[\tilde g(Z)
	\right]\right\vert\\ 
	&\leq	B_3\left(\tilde g,K_j\times d\right)\sum_{k=1}^{K_j}\sqrt{\operatorname{cum}_4(\hat \beta_{j,k;t})}a
	\end{split}
	\]
	where
	\[
	\begin{split}
		B_3(\tilde g,K_j\times d) & =B_3\left(g \circ R_{j,d},K_j\times d\right)\\
		&=\frac{\sqrt{2K_j}}{4}M_{2}(g\circ R_{j,d})+\frac{2\sqrt{K_j\operatorname{Tr}\Upsilon_d}}{9}M_3(g\circ R_{j,d})\\
		& \leq \frac{\sqrt{2}S_j}{4}M_{2}(g\circ R_{j,d})+\frac{2\bar \sigma}{9\sqrt{\varepsilon_j}}M_3(g\circ R_{j,d})
	\end{split}
	\]
	Now, observe that 
	\[
	\underset{\tilde g:\mathbb{R}^d\to \mathbb{R}}{ \sup}M_k (g \circ R_{j;d}) \leq \tilde{C}_M dS_j\varepsilon_j	\underset{g:\mathbb{R}^d\to \mathbb{R}}{ \sup}M_k (g ) ,
	\]
so that 
	\[
\begin{split}
\sup_{\tilde g}	B_3(\tilde g,K_j\times d) & \leq \sup_{g} \left(\tilde{C}_M ^\prime dS^2_j \varepsilon_j M_2(g) + C_M^{\prime \prime } d S_j \sqrt{\varepsilon_j} M_3(g)\right)\\
\leq O\left(dS^2_j \varepsilon_j \right)
	\end{split}
\]
So here we obtain a bound of order $O\left(\frac{dS_j ^3 \varepsilon_j}{\nu_t} \right)$, to be compared with the one previously obtained of order 
$O\left(\frac{d^2S_j ^2 \varepsilon_j}{\nu_t} \right)$.
Thus it holds that 
	\[
	\left \vert \mathbb{E} \left[g(\Xi_{;j,t})\right]-\mathbb{E}\left[g(Z_d)\right]\right \vert\leq O \left(\frac{dS_j ^2 \varepsilon_j \max(d,S_j)}{\nu_t}\right).	
	\]	
	The bound obtained here is consistent with the one presented in \cite{dmt24}[Remark 2.13]. Indeed, when considering an extremely shrinking needlet regime ($p>1$), the scaling $S_j$ effectively reduces to the classical multipole $\ell$ familiar from harmonic analysis, as discussed in \cite{d25}.
\end{remark}
\subsection{Functional convergence}
In the preceding sections, we focused on the limit behaviour of finite‑dimensional objects, specifically, vectors formed by Poisson needlet coefficients or evaluations at a fixed collection of $d$ points. A natural next step, inspired by \cite{bdmt24}, is to lift these results to the entire Poisson needlet field  $\{\Psi_{j;t}(\cdot)\}$ as a random function on $\mathbb{S}^2$. Crucially, since each $\Psi_{j;t}$ can be written as a weighted linear combination of the spherical Poisson waves $\{T_{\ell;t}\}$ analyzed in \cite{bdmt24}, we can directly leverage their functional limit results. Recent advances in the topic (see among others \cite{bc}) then enable us to derive functional convergence rates in both the Hilbert space $L^2(\mathbb{S}^2)$ and the Sobolev space $W_{\alpha,2}(\mathbb{S}^2)$ for any $\alpha>0$. We will measure distances in these spaces using the metrics $d_{3;L^2(\mathbb{S}^2)}$ and $d_{3;W_{\alpha,2}(\mathbb{S}^2)}.$
Recall also that for a general Banach space $K$, the class $C^3_b(K)$ consists of all real‑valued functions on $K$ whose first through third Fréchet derivatives are uniformly bounded, endowed with the norm
$$
\|h\|_{C^3_b(K)} 
= \underset{\iota=1,2,3}{\sup} \underset{x\in K}{\sup}\|D^\iota h(x)\|_{K^{\otimes j}}.
$$
\subsubsection*{Functional limit theorems in $L^2\left(\mathbb{S}^2\right)$}
Recall that here, for any $x \in \mathbb{S}^2$, we consider%
$$
\Psi_{j;t }(x)=\frac{1}{\sqrt{\nu_t}\sigma_j}\sum_{k=1}^{N_t\left(\mathbb{S}^2\right)}\sum_{\ell =S_{j-1}}^{S_{j+1}}b_j^2\left(\ell\right)\frac{2\ell +1}{4\pi }P_{\ell }(\left\langle x,\xi _{k}\right\rangle ),
$$
where $N_t\left(\mathbb{S}^2\right)$ is a Poisson process with intensity $\nu_t$. \begin{align*}
	\sigma_j^2 & 
	=\sum_{\ell =S_{j-1}}^{S_{j+1}}b_j^4\left(\ell\right) \frac{2\ell +1}{4\pi }.
\end{align*}
Note that the field is centered, i.e. for any $x \in \mathbb{S}^2$ it holds that $\mathbb{E}\left[\Psi_{j;t }(x)\right]=0$, and , for any $x,y \in \mathbb{S}^2$, its covariance function is given by
\begin{align}\notag
	\Gamma_{j;t}(x,y) &=\mathbb{E}\left[\Psi_{j;t }(x)\Psi_{j;t }(y)\right]\\ 
	&=\frac{1}{\sigma_j^2}  \sum_{\ell^\prime =S_{j-1}}^{S_{j+1}} b_j^4\left(\ell\right) Z_{\ell }(\left\langle x,y\right\rangle ),\label{eq:cov}
\end{align}
such that $\operatorname{Var}\left(\Psi_{j;t}(x)\right)=1$ for any $x \in \mathbb{S}^2$. 
\begin{theorem} \label{thm:func}
	Let $\mathcal{Z}$ be a centred Gaussian process with the same covariance operator as $\Psi_{j;t }(x)$. We have that
	$$
d_{3,L^2\left(\mathbb{S}^2\right)}  \left(\Psi_{j;t},Z\right)\leq \left(\frac{1}{4}+4\sqrt{\pi} \right) \sqrt{\frac{4\pi}{\nu_t}}.
	$$
\end{theorem}
\begin{remark}[Topological Effects on Gaussian Approximation Rates]
	As observed in \cite{bdmt24}, in this functional framework the convergence rate of order $\nu_t^{-\frac{1}{2}}$ 
does not depend on the resolution index $j$, and is therefore faster than in the finite‑dimensional setting. This mirrors the behavior of spherical Poisson waves, where the convergence rate is likewise independent of the multipole $\ell$. The reason is that the $L^2\left(\mathbb{S}^2\right)$
	topology is too coarse to enforce convergence of finite‑dimensional distributions—in the harmonic picture this manifests as an inability to control individual multipoles, and in the needlet decomposition it appears through weighted superpositions of those multipoles. In the next subsection, we turn to Sobolev spaces, which induce a finer topology and allow us to recover quantitative control over finite‑dimensional marginals.
\end{remark}
\subsubsection*{Functional limit theorems in $W_{\alpha;2}\left(\mathbb{S}^2\right)$}
We consider the random Poisson needlets taking values in Sobolev spaces $W_{\alpha;2}\left(\mathbb{S}^2\right)$,
with $\alpha > 0$, on the sphere, corresponding to the spaces of functions $f\in L^2\left(\mathbb{S}^2\right)$, with finite norm 
\[
\left\Vert f\right \Vert_{W_{\alpha;2}\left(\mathbb{S}^2\right)}= \sum_{\ell\geq0}\sum_{m=-\ell}^{\ell}\left(1 +\sqrt{\ell\left(\ell+1\right)} \right)^{2\alpha} \left \vert a_{\ell,m} \right \vert^2.
\]
Our main finding is the described by the next theorem.
\begin{theorem}\label{thm:sob}
	Let $Z$ be a centered Gaussian process with the same covariance operator as $\Psi_{j;t}$. Then it holds that
	\begin{equation*}
		d_{3;W_{2;\alpha}}	  \left(\Psi_{j;t},Z\right)\leq\frac{2\sqrt{\pi}}{\sigma_j^2 \sqrt{\nu_t}} \left(\sum_{\ell=S_{j-1}}^{S_{j+1}}b_j^4(\ell)\frac{2\ell+1}{4\pi}\left(1+\sqrt{\ell(\ell+1)}\right)^{2\alpha}\right).
	\end{equation*}
\end{theorem}
\begin{corollary}\label{cor:Sobolev}
	If $S_j$ is subexponential and regular (as in Equation \eqref{eqn:Sj}), thus 
	\begin{equation*}
		d_{3;W_{2;\alpha}}	  \left(\Psi_{j;t},Z\right)\leq	O\left(\frac{S_j^{2\alpha}}{\sqrt{\nu_t}}\right)
	\end{equation*}
\end{corollary}
	\begin{remark}[Linking Sobolev and regularity parameters]
	We now derive explicit conditions on the resolution index $j$ ensuring that 
	\[
	S_j^{2\alpha} = o\left( \sqrt{\nu_t} \right),
	\]
	which guarantees that $	d_{3;W_{2;\alpha}}	  \left(\Psi_{j;t},Z\right)$ converges to 0. This translates to the logarithmic condition
	\[
	\log S_j = o\left( \tfrac{1}{4\alpha} \log \nu_t \right).
	\]
	Recall that $\gamma(j)$ is a slowly varying function, i.e., asymptotically behaves like a constant up to logarithmic corrections. Then, the growth conditions on $j$ collected in Table \ref{tab:SobolevCondition} are sufficient in the respective regimes.
	
	\begin{table}[h!]
		\centering
		\renewcommand{\arraystretch}{1.3}
		\begin{tabular}{ll}
			\hline
			\textbf{Case} & \textbf{Condition on $j(t)$ for convergence} \\
			\hline
			Stretched exponential $(0 < p < 1)$ 
			& $j(t) = o\left( \left( \log \nu_t \right)^{\frac{1}{(1-p)\cdot 4\alpha}} \right)$ \\
			
			Logarithmic growth $(p = 1)$ 
			& $j(t) = \nu_t^{o(1)}$ \\
			\hline
		\end{tabular}
		\vspace{3pt}
		\caption{Conditions ensuring $S_j^{2\alpha} = o(\sqrt{\nu_t})$ under slowly varying $\gamma(j)$.}
		\label{tab:SobolevCondition}
	\end{table}
	These asymptotic regimes reflect the interplay between the scaling behavior of $S_j$ and the concentration rate of the Poisson intensity $\nu_t$, modulated by the Sobolev index $\alpha$. The required regularity on $j$ becomes more stringent as $\alpha$ increases. In the stretched exponential regime, the condition reflects the balance between the almost-exponential growth of $S_j$ and the regularity parameter $\alpha$, which penalizes high-frequency resolution. In contrast, the logarithmic case yields a more permissive criterion: the resolution level $j$ may grow almost polynomially with $\nu_t$, making it easier to satisfy in practice. In both regimes, increasing $\alpha$ (i.e., seeking convergence in a stronger Sobolev topology) forces $j$ to grow more slowly with $\nu_t$.
	\end{remark}
\section{Proofs}\label{sec:proofs}
\begin{proof}[Proposition \ref{prop:sigma}]
	First observe that, for any $n\in\mathbb{N}$ and $j \geq 1$, since $b_j$ is compactly supported, it holds that, for any $u \in(S_{j-1},S_{j+1})$,
	$$0<\kappa_n\leq b^n_j(u) \leq K_n<\infty$$. Hence it follows that
	\[
	\sum_{\ell \in [S_{j-1},S_{j+1}]} b_j^n(\ell) \frac{2\ell+1}{4\pi} \simeq \frac{C_{b,n}}{4\pi} \left(S_{j+1}^2-S_{j-1}^2\right)+o\left(S_{j+1}^2-S_{j-1}^2\right) 
	\]
	Then, note that 
	$$
	S_{j+1} =S_j \exp(\varepsilon_{j}); \quad	S_{j-1} =S_j \exp(\varepsilon_{j-1})
	$$
	Therefore, 
	\[
	\begin{split}
		S_{j+1}^2-S_{j-1}^2&=S_j^2\left( e^{2\varepsilon_j}-e^{2\varepsilon_{j-1}}\right)\\
		&= 2S_{j}^2\left( \left(\varepsilon_{j}+\varepsilon_{j-1}\right)+o\left(\varepsilon_{j}+\varepsilon_{j-1}\right) \right),
	\end{split}
	\]
	where in the last equality we use the Taylor expansion around 0, since $\varepsilon \rightarrow 0$.
	Now, recall that $\varepsilon_j= \gamma(j)/j^p$; since $\gamma$ is slowly varying we have that 
	$$
	\frac{\gamma(j-1)}{\gamma(j)} \rightarrow 1
	$$
	yielding
	\[
	\begin{split}
		\varepsilon _{j} + \varepsilon_{j-1}& = \frac{\gamma(j)}{j^p}+\frac{\gamma(j-1)}{(j-1)^p}\\
		&=	\gamma(j)\frac{1}{j^p}\left(1 + \frac{\gamma(j-1)}{\gamma(j)} \left(1-\frac{1}{j}\right)^{-p} \right) \\
		&=	2\gamma(j)\frac{1}{j^p}\left(1 +o(j)\right)\\
		&=	2\varepsilon_j\left(1 +o(j)\right),
	\end{split}
	\]	
	so that
		\[
		S_{j+1}^2-S_{j-1}^2=4S_j^2\varepsilon_j\left(1 +o(j)\right).
		\] 
		Therefore,
	\[
	\sigma^2_j \sim  \frac{C_b}{\pi} (S_j^2 \varepsilon_j),
	\]
	as claimed.

\end{proof}
	\begin{proof}[Proposition \ref{prop:dimone0}]
	Using \eqref{eqn:wass} mutuated from \cite{dvz} and resumed in Section \ref{sec:appuni} yields
	\[
	\begin{split}
		d_{\operatorname{W}} \left(\hat{\beta}_{j,k;t},Z \right)  & \leq \tilde{C}
		\sqrt{\operatorname{cum}_4\left(\hat{\beta}_{j,k;t}\right)} 
	\end{split},
	\]
	On the one hand,  $\operatorname{Var}\left(\hat{\beta}_{j,k;t}\right)=\frac{\bar{\sigma}^2}{S^2_j \varepsilon_j} $,
	while on the other 
	\[
	\begin{split}
		\mathbb{E}\left[\hat\beta_{j,k;t}^4\right] &\leq \frac{1}{S_j^4 \varepsilon_j^2}\frac{\bar{\sigma}^4}{\nu_t} \int_{\mathbb{S}^2} \left \vert \psi_{j,k}(x)\right \vert^4 \operatorname{d}x  + 3  \operatorname{Var}\left(\hat{\beta}_{j,k;t}\right)^2\\
		& \leq\frac{1}{S_j^4 \varepsilon_j^2}\frac{\bar{\sigma}^4}{\nu_t} \left \Vert \psi_{j,k}\right \Vert_4^4 + 3  \operatorname{Var}\left(\hat{\beta}_{j,k;t}\right)^2\\
		&  \leq\frac{1}{S_j^4 \varepsilon_j^2}\frac{\zeta_4\bar{\sigma}^4}{\nu_t}  \Sigma_{j;p}^{2}  + 3  \operatorname{Var}\left(\hat{\beta}_{j,k;t}\right)^2\\
		&   \leq\frac{1}{S_j^2 }\frac{\zeta^4_4\bar{\sigma}^4}{\nu_t}    + 3  \operatorname{Var}\left(\hat{\beta}_{j,k;t}\right)^2.
	\end{split}
	\]
	Combining these two findings yields the stated conclusion.
\end{proof}
\begin{proof}[Corollary \ref{cor:dimone0}]
	Using the results for \cite{dvz} resumed in Section \ref{sec:appuni} yields
	\[
	\begin{split}
		d_{\operatorname{W}} \left(\tilde{\beta}_{j,k;t},\tilde Z \right)  & \leq \tilde{C}
		\sqrt{\operatorname{cum}_4\left(\tilde{\beta}_{j,k;t}\right)} 
	\end{split},
	\]
	Recalling that $\operatorname{Var}\left(\tilde{\beta}_{j,k;t}\right)=1$, its fourth moment corresponds to
	\[
	\begin{split}
		\mathbb{E}\left[\tilde\beta_{j,k;t}^4\right] &\leq \frac{1}{\nu_t\tilde{\sigma}^4} \int_{\mathbb{S}^2} \left \vert \psi_{j,k}(x)\right \vert^4 \operatorname{d}x  + 3  \operatorname{Var}\left(\tilde{\beta}_{j,k;t}\right)^2\\
		& \leq \frac{1}{\nu_t\tilde{\sigma}^4}\left \Vert \psi_{j,k}\right \Vert_4^4 + 3\\
		&  \leq \frac{1}{\nu_t\tilde{\sigma}^4} \Sigma_{j;p}^{2}  + 3 \\
		&   \leq \frac{\zeta_4^4}{\tilde{\sigma}^4} \frac{S_j^2\varepsilon_j^2}{\nu_t} + 3. 
	\end{split}
	\]
	As above, a simple sostitution leads to the claimed results.
\end{proof}
\begin{proof}[Proposition \ref{prop:multicoef}]
	We use the results from \cite{dvz} recalled in Equation \eqref{thm:A3bis}, stressing that, here, $$\sum_{k=1}^{K_j}\sqrt{\mathbb{E}\left[\hat \beta_{i,k;t}^4\right]-3\mathbb{E}^2\left[\hat \beta_{i,k;t}^2\right]}=O\left(K_j\frac{1}{S_j \sqrt{\nu_t}}\right)=O\left(\frac{S_j}{\sqrt{\nu_t}}\right),$$ from Proposition \ref{prop:dimone0}. 
	Also, following Equation \eqref{eqn:Btre} yields
	$$
	B_3(g,S_j)=\frac{\sqrt{2K_j}}{4}M_2(g)+ \frac{2}{9}\sqrt{K_j \operatorname{Tr}\left(\Upsilon_{j;t}\right)}M_3(g).
	$$
	The claimed result follows by observing that by construction
	$\operatorname{Tr} \Upsilon_{j;t} = \bar{\sigma}^2/\varepsilon_j$, while $K_j = O(S_j^2)$. Finally since $\varepsilon_j$ tends to 0, to achieve convergence in law it is sufficient to ensure that $S_j^{4}/\nu_t\rightarrow 0$.
\end{proof}
	\begin{proof}[Corollary \ref{cor:multicoef}]
	We use again Equation \eqref{thm:A3bis} (see  \cite{dvz}). Observe that here it follows from Proposition \ref{prop:dimone0} that
	$$\sum_{k=1}^{K_j}\sqrt{\mathbb{E}\left[\tilde \beta_{i,k;t}^4\right]-3\mathbb{E}^2\left[\tilde \beta_{i,k;t}^2\right]}=O\left(K_j\sqrt{\frac{S^2_j \varepsilon^2_j}{\nu_t}}\right) = O\left(\sqrt{\frac{S^6_j \varepsilon^2_j}{\nu_t}}\right).$$
	Also, it follows from \eqref{eqn:Btre} that
	$$
	\tilde B_3(g,S_j)=\frac{\sqrt{2}S_j}{4}M_2(g) + \frac{2}{9}S_j^2M_3(g),$$
	since the covariance matrix $\tilde \Upsilon_{j;t} $ converges to the identity and then 
	$\operatorname{Tr}\tilde  \Upsilon_{j;t} = K_j$.
\end{proof}
\begin{proof}[Theorem \ref{thm:fin1d}]
	First of all, note that
	\begin{equation*}
		\begin{split}
			\mathbb{E}\left[ \Psi_{j;t}^4 \right] - 3 & = \frac{1}{\nu_t^2 \sigma_j^4} \mathbb{E}\left[ \int_{\left(\mathbb{S}^2\right)^{\otimes 4}} \prod_{i=1}^4 \Phi_j (x,y_i)N_t(\operatorname{d}y_i)\right] -3\\
			&=\frac{1}{\nu_t^2 \sigma_j^4} \left[ \nu_t \int_{\mathbb{S}^2}  \Phi^4_j (x,y) \operatorname{d}y + 3\nu_t^2 \left(\int_{\mathbb{S}^2}  \Phi^2_j (x,y) \operatorname{d}y  \right)^2\right]-3\\
			&\leq \frac{1}{\nu_t \sigma_j^4}  \int_{\mathbb{S}^2} \left \vert \Phi_j (x,y)\right \vert ^4 \operatorname{d}y \\
			&\leq \frac{\pi}{ \nu_t  S_j^4 \varepsilon_j^2  } 4\pi C \left(S_{j+1}^2-S_{j-1}^2\right)^4 		
		\end{split}
	\end{equation*}
	where in the last inequality we have used the localization property \eqref{eqn:locPhi}.  Now, recall that  
	\begin{equation*}
		\begin{split}
			S_{j+1}^2-S_{j-1}^2 & =S_j^2\left( \exp \left( 2 \varepsilon_{j}\right)-\exp \left( 2 \varepsilon_{j-1}\right)\right)\\
			&= 2 S_j^2\left( \varepsilon_j + \varepsilon_{j-1} \right)\left(1+ o(j)\right)\\
			&= 4 S_j^2 \varepsilon_j  (1+ o(j)).
		\end{split}
	\end{equation*}
	Thus it follows that
	\begin{equation*}
		\begin{split}
			\mathbb{E}\left[ \Psi_{j;t}^4 \right] - 3 &\leq 4^{5}\pi^2\frac{S_j^4 \varepsilon_j^2 }{ \nu_t   } (1+ o(j)),
		\end{split}
	\end{equation*}
	as claimed.
\end{proof}
\begin{proof}[Theorem \ref{thm:multi}]
	Since each $\Psi_{j}(x_i)$ lies in the first Poisson Wiener chaos, Equation \eqref{thm:A3} gives
	\[
	d_3(\Xi_{d;j,t} ,Z_d) \leq 
	\sup_{g\in C^3}B_3(g;d)\sum_{i=1}^d\sqrt{\mathrm{cum}_4\bigl(F_{d,i}\bigr)}.
	\]
	By Theorem \ref{thm:fin1d} each fourth cumulant satisfies
	\[\mathrm{cum}_4(\Psi_{j}(x_i))=\sqrt{\frac{S_{j}^4 \varepsilon^2_j}{\nu_t}}+ o\left(\sqrt{\frac{S_{j}^4 \varepsilon^2_j}{\nu_t}}\right).\] Finally, summing over $i=1,\dots,d$ yields the stated bound.  As far as $B_3(g,d)$ is concerned, it suffices to use Equation \eqref{eqn:Btre} whereas the covariance matrix is the identity of dimension $d$.
\end{proof}
\begin{proof}[Theorem \ref{thm:func}]
	Using \eqref{eq:cov}, the Hilbert-Schmidt norm of the covariance operator $S_{j;t}$ is given by
	\begin{align}\notag
		\left\Vert S_{j;t}\right\Vert_{\operatorname{HS}(L^2)}^2 & =\left \Vert \mathbb{E}\left[\Psi_{j;t}\otimes \Psi_{j;t} \right]     \right \Vert_{\operatorname{HS}}^2=\left \Vert \Gamma_{j;t}\right \Vert_2^2 \\ &=	\int_{\mathbb{S}^2}\int_{\mathbb{S}^2}\left \vert \frac{1}{\sigma_j^2}  \sum_{\ell=S_{j-1}}^{S_{j+1}} b_j^4\left(\ell\right) Z_{\ell }(\left\langle x,y\right\rangle ) \right \vert^2 \operatorname{d}x \operatorname{d}y \notag\\
		&=\frac{4\pi}{\sigma_j^4}  \sum_{\ell=S_{j-1}}^{S_{j+1}} b^8\left(\ell\right)\frac{2\ell +1}{4\pi }\label{eq:HSCov}
	\end{align}
First, note that
	\begin{align*}
		\mathbb{E}\left[\left\Vert\Psi_{j;t}\right \Vert_{L^2}^2\right] & = \mathbb{E}\left[\int_{\mathbb{S}^2}\left \vert \Psi_{j;t}\left(x\right)\right \vert ^2 \operatorname{d}x \right]\\
		& =\frac{1}{\nu_t\sigma_j^2} \mathbb{E}\left[\sum_{\iota,\iota^\prime=1}^{N_t\left(\mathbb{S}^2\right)}\sum_{\ell,\ell^\prime=S_{j-1}}^{S_{j+1}} 
		b_j^2\left(\ell\right)b_j^2\left(\ell^\prime\right) \int_{\mathbb{S}^2}Z_{\ell }(\left\langle x,\upsilon _{\iota}\right\rangle )Z_{\ell^\prime }(\left\langle x,\upsilon _{\iota^\prime}\right\rangle ) \operatorname{d}x \right]\\
		&=\frac{1}{\sigma_j^2} \sum_{\ell=S_{j-1}}^{S_{j+1}}b_j^4\left(\ell\right)\frac{2\ell +1}{4\pi }\int_{\mathbb{S}^2} P_{\ell }(\left\langle z ,z\right\rangle ) \operatorname{d}z \\ 
		&+ \frac{1}{\sigma_j^2} \sum_{\ell=S_{j-1}}^{S_{j+1}}b_j^4\left(\ell\right)\frac{2\ell +1}{4\pi }\int_{\mathbb{S}^2}\int_{\mathbb{S}^2} P_{\ell }(\left\langle z_1 ,z_2\right\rangle ) \operatorname{d}z_1	\operatorname{d}z_2\\
		&=\frac{4\pi}{\sigma_j^2} \sum_{\ell=S_{j-1}}^{S_{j+1}}b_j^4\left(\frac{\ell}{B^j}\right)\frac{2\ell +1}{4\pi }\\
		&=4\pi.	
	\end{align*}
Then, we have that
	\begin{align*}
		\mathbb{E}\left[\left\Vert\Psi_{j;t}\right \Vert_{L^2}^4\right]&=\mathbb{E}\left[\left\Vert\Psi_{j;t}\right \Vert_{L^2}^2\left\Vert\Psi_{j;t}\right \Vert_{L^2}^2\right] \\
		& = \mathbb{E}\left[\int_{\mathbb{S}^2}\left \vert \Psi_{j;t}\left(x\right)\right \vert ^2 dx \int_{\mathbb{S}^2}\left \vert \Psi_{j;t}\left(y\right)\right \vert ^2 \operatorname{d}y\right]\\
		& =\frac{1}{\nu_t^2\sigma_j^4} \mathbb{E}\left[\sum_{\iota_1,\iota_2,\iota_3,\iota_4=1}^{N_t\left(\mathbb{S}^2\right)}\sum_{\ell_1,\ell_2,\ell_3,\ell_4=S_{j-1}}^{S_{j+1}} \prod_{i=1}^4b_j^2\left(\ell\right)\right.\\&\left. \times \int_{\mathbb{S}^2}Z_{\ell_1 }(\left\langle x,\upsilon _{\iota_1}\right\rangle )Z_{\ell_2 }(\left\langle x,\upsilon _{\iota_2}\right\rangle ) \operatorname{d}x \int_{\mathbb{S}^2}Z_{\ell }(\left\langle y,\upsilon _{\iota_3}\right\rangle )Z_{\ell_4 }(\left\langle y,\upsilon _{\iota_4}\right\rangle ) \operatorname{d}y\right]
	\end{align*}
	Using the reproducing kernel property of $Z_\ell$, it holds that
	\begin{align*}
		\mathbb{E}\left[\left\Vert\Psi_{j;t}\right \Vert_{L^2}^4\right]
		& =\frac{1}{\nu_t^2\sigma_j^4} \mathbb{E}\left[\sum_{\iota_1,\iota_2,\iota_3,\iota_4=1}^{N_t\left(\mathbb{S}^2\right)}\sum_{\ell,\ell^\prime=S_{j-1}}^{S_{j+1}}b_j^4\left(\ell\right)b_j^4\left(\ell^\prime\right) Z_{\ell}(\left\langle \upsilon _{\iota_1},\upsilon _{\iota_2}\right\rangle ) Z_{\ell^\prime }(\left\langle \upsilon _{\iota_3},\upsilon _{\iota_4}\right\rangle )\right]
		\end{align*}
		Computing the expectation we obtain then three terms, such that
		\begin{align*}
			\mathbb{E}\left[\left\Vert\Psi_{j;t}\right \Vert_{L^2}^4\right]&= \frac{1}{\nu_t\sigma_j^4} \sum_{\ell,\ell^\prime=S_{j-1}}^{S_{j+1}}b_j^4\left(\ell\right)b_j^4\left(\ell^\prime\right)\int_{\mathbb{S}^2}Z_{\ell}(\left\langle z,z\right\rangle ) Z_{\ell^\prime }(\left\langle z,z\right\rangle )\operatorname{d}z\\&+ \frac{1}{\sigma_j^4} \sum_{\ell,\ell^\prime=S_{j-1}}^{S_{j+1}}b_j^4\left(\ell\right)b_j^4\left(\ell^\prime\right)\int_{\mathbb{S}^2}Z_{\ell}(\left\langle z,z\right\rangle )\operatorname{d}z \int_{\mathbb{S}^2} Z_{\ell^\prime }(\left\langle z^\prime,z^\prime\right\rangle )\operatorname{d}z^\prime\\
		&+ \frac{2}{\sigma_j^4} \sum_{\ell,\ell^\prime=S_{j-1}}^{S_{j+1}}b_j^4\left(\ell\right)b_j^4\left(\ell^\prime\right)\int_{\mathbb{S}^2}\int_{\mathbb{S}^2}Z_{\ell}(\left\langle z,z^\prime\right\rangle )Z_{\ell^\prime}(\left\langle z,z^\prime\right\rangle )\operatorname{d}z\operatorname{d}z^\prime \\
		&= \frac{4\pi}{\nu_t\sigma_j^4}\left( \sum_{\ell=S_{j-1}}^{S_{j+1}}b_j^4\left(\ell\right)\frac{2\ell+1}{4\pi}\right)^2
		+ \frac{(4\pi)^2}{\sigma_j^4} \left( \sum_{\ell={S_{j-1}}}^{S_{j+1}}b_j^4\left(\ell\right)\frac{2\ell+1}{4\pi}\right)^2\\
		&+ \frac{2}{\sigma_j^4} \sum_{\ell=S_{j-1}}^{S_{j+1}}b_j^8\left(\ell\right)\frac{2\ell+1}{4\pi}\int_{\mathbb{S}^2}P_{\ell}(\left\langle z,z\right\rangle )\operatorname{d}z\\
		&= \frac{4\pi}{\nu_t}+(4 \pi)^2 + \frac{8\pi}{\sigma^4_j}\sum_{\ell=S_{j-1}}^{S_{j+1}}b^8\left(\frac{\ell}{B^j}\right)\frac{2\ell+1}{4\pi}
	\end{align*}
The second term exactly offsets the square of the second‐order moment, whereas the third term vanishes against the Hilbert–Schmidt norm of $S_{j;t}$.
	Then, we have that
	$$
	d_{3,L^2\left(\mathbb{S}^2\right)} (\Psi_{j;t},Z)\leq \left(\frac{1}{4}+4\sqrt{\pi} \right) \sqrt{\frac{4\pi}{\nu_t}}.
	$$
\end{proof}
\begin{proof}[Theorem \ref{thm:sob}]
	Observe preliminarily that
	\[
	\begin{split}
		\Psi_{j;t}(x) &= \frac{1}{\sigma_j} \sum_{\ell=S_{j-1}}^{S_{j+1}} b_j^2(\ell)  \sqrt{\frac{2\ell+1}{4\pi}} T_{\ell;t}(x)\\
		&= \sum_{\ell=S_{j-1}}^{S_{j+1}} \omega_{j}(\ell)T_{\ell;t}(x).
	\end{split} 
	\]
	Thus, it follows that 
	\[
	\begin{split}
		\left\Vert \Psi_{j;t}(x)\right \Vert_{W_{\alpha;2}\left(\mathbb{S}^2\right)}= \sum_{\ell=S_{j-1}}^{S_{j+1}} \omega^2_j(\ell) \left\Vert T_{\ell;t}(x)\right \Vert_{W_{\alpha;2}}^2\\
	\end{split}
	\]
	First, since $	\mathbb{E}\left[ \left\Vert T_{\ell;t} \right \Vert_{W_{\alpha;2}}^2\right]=4\pi\left(1+\sqrt{\ell(\ell+1)}\right)^{2\alpha}$ (see \cite{bdmt24}), we have on the one hand that 
	\[
	\mathbb{E}\left[ \left\Vert \Psi_{j;t} \right \Vert_{W_{\alpha;2}}^2\right] =\frac{4\pi}{\sigma_j^2}\sum_{\ell=S_{j-1}}^{S_{j+1}} b_j^4(\ell) \frac{2\ell+1}{4\pi} \left(1+\sqrt{\ell(\ell+1)}\right)^{2\alpha},
	\]
	while on the other one
	\[
	\begin{split}
		\mathbb{E}\left[ \left\Vert \Psi_{j;t} \right \Vert_{W_{\alpha;2}}^4\right] &=	\mathbb{E}\left[ \left\Vert \Psi_{j;t} \right \Vert_{W_{\alpha;2}}^2\left\Vert \Psi_{j;t} \right \Vert_{W_{\alpha;2}}^2\right]\\
		&  =\sum_{\ell=S_{j-1}}^{S_{j+1}}\sum_{\ell^\prime=S_{j-1}}^{S_{j+1}} \omega_j^2(\ell) \omega_j^2(\ell^\prime)  \mathbb{E}\left[ \left\Vert T_{\ell;t} \right \Vert_{W_{\alpha;2}}^2\left\Vert T_{\ell^\prime;t} \right \Vert_{W_{\alpha;2}}^2\right],
	\end{split}
	\]
	Straightforward calculations show that 
	\[
	\begin{split}
		\mathbb{E}\left[ \left\Vert T_{\ell;t} \right \Vert_{W_{\alpha;2}}^2\left\Vert T_{\ell^\prime;t} \right \Vert_{W_{\alpha;2}}^2\right] &= \left(1+\sqrt{\ell(\ell+1)}\right)^{2\alpha}\left(1+\sqrt{\ell^\prime(\ell^\prime+1)}\right)^{2\alpha}\\
		&\times \frac{4\pi}{\nu_t} + (4\pi)^2 + 2\cdot4\pi \cdot \frac{4\pi}{2\ell+1}\delta_{\ell}^{\ell^\prime}.
	\end{split}
	\]
	Thus we obtain
	\[
	\begin{split}
		\mathbb{E}\left[ \left\Vert \Psi_{j;t} \right \Vert_{W_{\alpha;2}}^4\right] &=
		\frac{1}{\sigma_j^4}\left(\frac{4\pi}{\nu_t}+(4\pi)^2\right) \left(\sum_{\ell=S_{j-1}}^{S_{j+1}}b_j^4(\ell)\frac{2\ell+1}{4\pi}\left(1+\sqrt{\ell(\ell+1)}\right)^{2\alpha}\right)^2\\		
		&+ 2\cdot 4\pi \sum_{\ell=S_{j-1}}^{S_{j+1}}b_j^8(\ell)\frac{2\ell+1}{4\pi}\left(1+\sqrt{\ell(\ell+1)}\right)^{4\alpha}
	\end{split}
	\]
	Finally, using the results for the Hilbert Schmidt norm in the Sobolev case of covariance functions related to Spherical Poisson wave, the relation between them and the spherical Poisson needlets and Equation \eqref{eq:HSCov}, it follows that
	\[
	\left\Vert S_{j;t}\right\Vert^2_{\operatorname{HS}(W_{\alpha;2})} = 4\pi \sum_{\ell=S_{j-1}}^{S_{j+1}} b^8_j(\ell)\frac{2\ell+1}{4\pi} \left(1+\sqrt{\ell(\ell+1)}\right)^{2\alpha}.
	\]
	After some easy simplification, we have that
	\begin{equation*}
		d_{3;W_{2;\alpha}}	  \left(\Psi_{j;t},Z\right)\leq	\frac{2\sqrt{\pi}}{\sigma_j^2 \sqrt{\nu_t}} \left(\sum_{\ell=S_{j-1}}^{S_{j+1}}b_j^4(\ell)\frac{2\ell+1}{4\pi}\left(1+\sqrt{\ell(\ell+1)}\right)^{2\alpha}\right),
	\end{equation*}
	as claimed.
\end{proof}
\begin{proof}[Corollary\ref{cor:Sobolev}]
	First of all, observe that 
	\[
	\sigma_{j}^4 \sim S_j^4\varepsilon_j^2.
	\]
	Also observe that 
	\[
	1+\sqrt{\ell(\ell+1)} \approx \sqrt{\ell^2+\ell}\approx \ell+\frac{1}{2}+o(1), 
	\]
	so that
	\[
	(1+\sqrt{\ell(\ell+1)} )^{2\alpha} \sim \ell^{2\alpha}, 
	\]
	while $2\ell+1 \approx 2\ell$. We can then easily compute that
	\[
	\sum_{\ell=S_{j-1}}^{S_{j+1}}b_j^4(\ell)\frac{2\ell+1}{4\pi}\left(1+\sqrt{\ell(\ell+1)}\right)^{2\alpha}=\frac{1}{2\pi(2\alpha+2)} (S_{j+1}^{2\alpha+2}-S_{j-1}^{2\alpha+2})(1+o(1)).
	\]
	Now, standard facts used in the proof of Proposition \ref{prop:sigma} yields
	\[
	(S_{j+1}^{2\alpha+2}-S_{j-1}^{2\alpha+2})=S_{j}^{2\alpha+2}2(2\alpha+2)\varepsilon_j. 
	\]
	Straightforward calculations lead to the claimed result.
\end{proof}
	\section*{Funding}
This work was partially supported by the PRIN 2022 project GRAFIA (Geometry of Random Fields and its Applications), funded by the Italian Ministry of University and Research (MUR).

\end{document}